\documentclass[12pt]{article}
\usepackage{amsfonts,amssymb,latexsym}
\normalbaselines
\begin{document}

\newcommand{\PP}{{\cal P}}
\newcommand{\XX}{{\cal X}}
\newcommand{\OO}{{\cal O}}
\newcommand{\DD}{{\cal D}}
\newcommand{\FF}{{\cal F}}
\newcommand{\II}{{\cal I}}
\newcommand{\JJ}{{\cal J}}
\newcommand{\QQ}{{\cal Q}}
\newcommand{\MM}{{\cal M}}
\newcommand{\Cb}{\Bbb C}
\newcommand{\Sb}{{\Bbb S}}
\newcommand{\Qb}{{\Bbb Q}}
\newcommand{\eps}{\varepsilon}
\newcommand{\Zb}{\Bbb Z}
\newcommand{\Tb}{\Bbb T}
\newcommand{\Nb}{\Bbb N}
\newcommand{\Eb}{\Bbb E}
\newcommand{\Mb}{\Bbb M}
\newcommand{\Fract}{{\rm{Fract}}}
\newcommand{\Ker}{{\rm{Ker}}}
\newcommand{\Spec}{{\rm{Spec}}}
\newcommand{\Ad}{{\rm{Ad}}}
\newcommand{\rank}{{\rm{rank}}}
\newcommand{\mod}{{\rm{mod}}}
\newcommand{\gr}{{\rm{gr}}}
\newcommand{\Maxspec}{{\rm{Maxspec}}}
\newcommand{\um}{\underline{m}}
\newcommand{\ual}{\underline{\alpha}}
\newcommand{\un}{\underline{n}}
\newcommand{\Rchi}{R_{\eps,\chi}}
\newcommand{\tx}{\tilde{x}}
\newcommand{\tm}{\tilde{m}}
\newcommand{\wa}{\widehat{a}}
\newcommand{\wT}{\widehat{T}}
\newcommand{\wx}{\widehat{x}}
\newcommand{\wR}{\widehat{R}}
\newcommand{\tR}{\tilde{R}}
\newcommand{\tRe}{\tilde{R}_\eps}
\newcommand{\ta}{\tilde{a}}
\newcommand{\tr}{\tilde{r}}
\newcommand{\ts}{\tilde{s}}
\newcommand{\tSb}{\tilde{\Sb}}
\newcommand{\tM}{{\tilde{M}}}
\newcommand{\tA}{{\tilde{A}}}
\newcommand{\tID}{{\tilde{I}(\pi,\DD)}}
\newcommand{\tI}{{\tilde{I}}}
\newcommand{\ttau}{{\tilde{\tau}}}
\newcommand{\tdelta}{{\tilde{\delta}}}
\newcommand{\tmD}{{\tilde{m}(\pi,\DD)}}
\newcommand{\tZe}{{\tilde{Z}_\eps}}
\newcommand{\tZ}{{\tilde{Z}}}
\newcommand{\pg}{{\frak p}}
\newcommand{\nog}{{\frak n}}
\newcommand{\tg}{{\frak t}}
\newcommand{\Jg}{{\frak J}}
\newcommand{\vg}{{\frak v}}
\newcommand{\kg}{{\frak k}}
\newcommand{\ug}{{\frak u}}
\newcommand{\Ng}{{\frak N}}
\newcommand{\qg}{{\frak q}}
\newcommand{\bg}{{\frak b}}
\newcommand{\Fg}{{\frak F}}
\newcommand{\Lg}{{\frak L}}
\newcommand{\mg}{{\frak m}}
\newcommand{\FAq}{{\mbox{FA}}_q}
\newcommand{\ua}{{\underline{a}}}
\newcommand{\DDJ}{{_\DD J}}
\newcommand{\Hle}{K[h_{1\eps}^{\pm l},\ldots,
h_{2t,\eps}^{\pm l}]}
\newcommand{\ZZ}{{\cal Z}}
\newcommand{\jg}{{\frak j}}
\newcommand{\Yg}{{\frak Y}}
\newcommand{\Zg}{{\frak Z}}
\newcommand{\gog}{{\frak g}}
\newcommand{\Sg}{{\frak S}}
\newcommand{\Rg}{{\frak R}}
\newcommand{\Ig}{{\frak I}}
\newcommand{\YY}{{\cal Y}}
\newcommand{\RR}{{\cal R}}
\newcommand{\tZY}{{\tZ_\eps\bigcap Y_\eps}}
\author{A.~N.~Panov.}
\date{\it{Mathematical Department, Samara State University,\\
 ul.Akad.Pavlova 1, Samara, 443011, Russia\\
panov@ssu.samara.ru}}
\title{
Irreducible representations of quantum  solvable algebras at roots of 1}
\maketitle
\footnote{The work is supported by RFFI grant  02-01-00017.}

\begin{abstract}

We study the irreducible representations of quantum solvable algebras
at roots of 1 which lie over a
point of the variety of center. 
We characterize the quiver of fiber algebra
and present the formulas on the dimension and 
the number of these representations in terms of Poisson 
structure of the variety of center.

\end{abstract}

\section {Introduction.}
Quantum algebras appears in papers on mathamatical physics as
deformations of the algebra
of regular functions
$\Cb[G]$ on the Lie group and universal enveloping algebra $U(\gog)$.
From algebraic point of view, quantizing $\Cb$-algebra $R$, we have got 
$\Cb$-algebra $R_q$ which is a free module over the ring of Laurent polynomials
$\Cb[q,q^{-1}]$ and $R = R_q\bmod(q-1)$.
If $R$ is a Hopf algebra, then it is natural to seek its quantizations 
in the class
of Hopf algebras.
The most familiar quantum algebras are 
quantum universal enveloping algebra $U_q(\frak{g})$ 
for semisimple Lie algebra $\frak{g}$, 
its dual Hopf algebra $\Cb_q[G]$, algebra of Quantum matrices, 
Quantum Weyl algebra.
 One can extend the chain of examples 
 considering the multiparatmeter versions of these algebras, 
 quantum spaces of representations.  
 
One sets up the problem of description 
of the space of primitive ideals.
It is itreresting to construct some general theory in spirit 
of the orbit method and also to classify primitive ideals for 
specific quantum algebras.
The problem reduces to specializations
 $R_\eps=R_q\bmod(q-\eps)$ where $\eps\in\Cb$.
  Two cases take place: $\eps$ is a root of 1  and $\eps$ is not a root of 1.
  
Up today the classification 
of primitive ideals is
known for for $\Cb_q[G]$ and quantum universal enveloping algebra 
of maximal solvable (resp. nilpotent) subalgebra in  $\frak{g}$.
The case of not a root of unity 
is studied in the book [J].
The papers  [DC-K], [DCKP1,2], [DC-L], [DC-P1,2]
are devoted to the case $\eps$  is a root of  1.

Here is the simplified plan of classification
of primitive ideals for  $\Cb_q[G]$.
The classification is based on  the description 
of symplectic leaves  on $G$ as orbits of dressing
transformations [ST].
For any symplectic leaf $\Omega$ one considers the ideal of functions 
vanishing on it. 
Its generators are some matrix elements of irreducible
representations of the  Lie group $G$. 
One can construct quantum analog of this ideal as an ideal 
generated by the corresponding matrix elements
of irreducible representations of  $U_q(\frak{g})$.
  The constructed ideal is primitive, if  $\eps$ is not a root of 1.
  It helps to stratify primitive ideals, if
 $\eps$ is a root of 1.
 
 The next example is algebra of Quantum matrices.
 Theses algebra is a bialgebra, but not a Hopf algebra.
 the above methods are not valid for it. For classification of 
 prime winding-invariant ideals see [GLn1,2], [C2], [L].
 One of the main goals is $U_q(\frak{g})$.
 This problem is far from its solution [J], [DC-K].

Consideration of examples make possible to 
set up some conjectures.
The next goal is to prove this conjectures  in 
maximally  weak assumptions imposed  on $R_q$. 
These assumptions must be easily checkable
and the theory must cover the main examples.
  
This paper is devoted to the case  of roots 1.
In what follows we suppose that  $\eps$ is a primitive  $l$th root of 1.
In the above examples, $R_\eps$ is finite over its center. 
That is the algebra $R_\eps$ is an order. 
Notice that  this property also holds for
elliptic algebras [FO],
some new quantum groups that appear in the framework
of theory of special functions [IK], 
reflection algebas [BG1].

The problem of description of primitive ideals 
 for orders  is equivalent to promlem of  
 classification of irreducible representations.
The restriction on the center of an irreducible representation
 $\pi$ of  $R_\eps$
is scalar  
$\pi\vert_{Z_\eps}=\chi\cdot\rm{id}$ 
and it defines the character  $\chi$
(i.e. the point of the variety)  of center $Z_\eps$.
We set up the usual problem for orders:
to classify all irreducible representations
of $R_\eps$ lying over given point $\chi$ of the variety of center.

There is one common feature of the above orders: 
the existence of the quantum adjoint action (see Section 2 and [DCKP2],[P3]).
Acting on the center $Z_\eps$, the quantum adjoint 
action defines the Poisson bracket.
The variety of center becomes a Poisson variety which splits 
into symplectic leaves.
It is proposed that the problem of classificationof irreducible representations
 can be solved in terms   
of geometrical and Poisson properties of the variety of center.

In the paper we study the quantum solvable algebras which are 
iterated skew polynomials extensions of $K[q,q^{-1}]$.
The examples of these
algebras are the algebra of Quantum matrices (see 2.14),
Quantum Weyl algebra, $U_q(\frak{b})$ and $U_q(\frak{n})$ (see.2.15) and their
numerous subalgebras.
The algebra  $\Cb_q[G]$ is not solvable, but 
one can reduce it to some solvable algebra after the localizaton.
For details in examples see [P2].
The main goal is the construction of quantum version 
of theory of Dixmier for $U_q(\frak{g})$ where $\frak{g}$ is a solvable 
Lie algebr [D].   
Here are some problems which stimulate general theory.\\
Problem 1. To prove that the symplectic leaves are algebraic
(i.e  Zariski-open in its Zariski closure); \\
Problem 2. To prove that the dimension of
an irreducible representation over $\chi$ is equal to  $l^{\frac{d}{2}}$
where $d$ is the dimension of symplectic leaf of $\chi$.
Conjectured in [DCP1,4.5],[DCP 25.1];\\
Problem 3. 
To describe the quiver of the algebra
$R_{\eps,\chi}:=R_\eps/m(\chi)R_\eps$ where $m(\chi)=\Ker(\chi)$.\\
Problem 4. To find the formula for 
the number of irreducible representations
over $\chi$.  

The solutions is known for  $\Cb_q[G]$ and $U_q(\frak{b})$.
The solution of Problem 1 for these algebras arises from 
the method of dressing transformations.
The formulas on dimensions and the number of irreducible 
representations  were obtained  in  [DC-P2]. The quivers were 
studed in [BG2].

In [P3] the Problems 1 and 2 were solved for rather
great $l$ (the point of good reduction of stratification process).
The goal of this paper is to drop these undesirable restriction on $l$
and to go forward in describing the quiver and determining the number of 
irreducible representations
over $\chi$.

The main definition of the paper is the definition of 
normal quantum solvable algebra (or NQS-algebra, see Definition 2.10).
We require that this algebra obeys some Conditions CN1,CN2.
We present two examples (Quantum matrices and $U(\nog)$).
One can find the other examples in ([G],[P1-P2]).
Our definition of admissible $l$  (Definition 2.18) is easily checkable
and necessary for solution of Promlems 1-4.

We stratify the prime $\DD$-stable spectrum
of NQS-algebra (see Theorem 3.2).
It is proved that every
prime $\DD$-stable ideal is completely prime (see Theorem 3.3).
The Problems 1 and 2 are solved in Theorem 4.2.

One can correspond the quiver to any finite dimensional 
algebra $A$ [Pie, 6.4].
 The vertices of quiver are primitive
 idempotents $e_1,\ldots,e_N$ such that their
 right ideals  $e_1A,\ldots, e_NA$ 
   represent non isomorphic classes of principal indecomposable
    $A$-modulas. 
 Two vertices $e_i, e_j$ are  linked with wedge $(e_i,e_j)$ if
 $e_iJe_j\ne 0$ where $J$ is the radical of $A$.
 In the paper we prove (see Theorem 4.3) that 
  any two vertices $e_i, e_j$ of quiver of 
  finite dimensional  algebra $R_{\eps,\chi}$
 is  linked by  wedges   $(e_i,e_j)$ and $(e_j,e_i)$. 
 In particular the quiver is connected.

In the last Section 5, we prove (Theorems 5.5,5.7) 
that the number of irreducible
representations over $\chi$ is 
 equal to $l^t$ where $t$ is the dimension of some toric Lie subalgebra 
of the stabilizer $\gog(\chi)$ of $\chi$ (Definition 5.6).

We are very thankfull to J.Cauchon; he sent his new 
preprint [C1] to the author.
The method of stratification of [C1] is used in this paper.
We are very thankfull to C.De Concini, C.Procesi, K.Brown and I.Gordon 
for useful discussions.

\section{Quantum solvable algebras and FA-elements}
We begin with some general definitions and the properties of skew extensions 
which are used throughout this paper.

Let $R_F$ be a domain and an algebra over
a field $F$.\\
{\bf Definition 2.1.} We say that $x\in R_F$ is an
element of finite adjoint action  (or $x$ is a FA-element) if 
$x$ is not a zero divisor and for every $a\in R_F$ there exists a polynomial
$f_a(t)= c_0t^N + c_1t^{N-1} +\cdots + c_N$, $c_0\ne 0$, $c_N\ne 0$ 
over $F$ such that 
$$ 
c_0x^Na + c_1x^{N-1}ax +\cdots + c_Nax^N = 0.\eqno (2.1) 
$$ 
A FA-element $x$ generates a denominator set $S_x:=\{x^n\}_{n\in\Nb}$
[P1, Proposition 3.3]. One can rewrite (2.1) in the form
$$
f_a(\Ad_x)a=0\eqno(2.2)
$$
where $\Ad_x(a)=xax^{-1}$. If $x$ is a FA-element in $R$, then 
it is a FA-element in $RS_x^{-1}$. The following statements are easy
to prove.\\
{\bf Proposition 2.2}.
Let $x,y\in\Fract(R_F)$ be FA-elements in a domain $R_F$ and suppose
that $xy=\gamma yx$ with some $\gamma\in F^*$. Then
$xy$ is also a FA-element.\\
{\bf Proposition 2.3}. 
Suppose that the above domain $R_F$ is
generated by $x_1,\ldots, x_n$ and $x\in \Fract(R_F)$.
Suppose that for every $j$ there exists a 
polynonial $f_j(t)$ obeying (2.1) with
 $a=x_j$. Then $x$ is a FA-element in $R_F$.
If, in addition,
$f_j(t)$ splits $f_j(t)=(t-\gamma_j^{(1)})\cdots (t-\gamma_j^{(n_j)})$,
then, for any $a\in R_F$, the polynomial $f_a(t)$ also splits with
the roots
in the semigroup generated by $\gamma_j^{(s)}$.

Let we have an endomorphism $\tau$ of $R_F$ 
($\tau$ is identical on $F$) and a $\tau$-derivation
$\delta$ of $R_F$ (i.e. $\delta(ab)=\delta(a)b+\tau(a)\delta(b)$ for all 
$a,b\in R_F$)
which is zero on $F$.
An Ore extension (skew extention)
$T_F=R_F[x;\tau,\delta]$ of $R_F$ is generated by $x$ and $R_F$ with
$xa=\tau(a)x+\delta(a)$ for all $a\in R_F$.
Every element of $T$ can be uniquely presented in the form
$\sum x^ir_i$ ( or $\sum r_ix^i$) where $r_i\in R$.
\\
{\bf Proposition 2.4}. Let $R_F$ and $T_F=R_F[x;\tau,\delta]$ be 
as above with
diagonalizable automorphism $\tau$. Suppose that
$\tau\delta=\gamma\delta\tau$, $\gamma\ne 0$.
The element $x$ is a FA-element in $T_F$ iff $\delta$ is locally nilpotent.
Moreover, for $\tau$-eigenvector $a$, there exists a polynomial 
$f_a(t)$ of degree $N$ obeying (2.1) iff $\delta^N(a)=0$.\\
{\bf Proof}. Let $a$ be a $\tau$-eigenvector, i.e. $\tau(a)=\beta a$.  
There exists a polynomial $f(t)$ obeying (2.1).
Then 
$$0=c_0x^Na+c_1x^{N-1}ax+\ldots+c_Nax^N=f(\beta) ax^N+\{\mbox{ 
 terms of lower
degree}\}.$$
It implies that $f(\beta)=0$,  
$f(t)=f_1(t)(t-\beta)$ and 
$0= f(\Ad_x)a= f_1(\Ad_x)(\Ad_x-\beta)a= f_1(\Ad_x)\delta(a)x^{-1}.$
The element $\delta (a)$ is also a $\tau$-eigenvector.
After $N$ steps we get $\delta^N(a)=0$
where $N=\deg f(t)$. 
On the other hand, if $\delta^N(a)=0$ and $\tau(a)=\beta a$, then
the polynomial 
$$
f(t):=\prod_{i=1}^N(t-\beta\gamma^i)
$$
obeys (2.1). 
$\Box$\\
Let $K$ be an algebraic closed field of zero characteristic,
 $q$ be an indeterminate and $C$ be a localization $K[q,q^{-1}]$ 
over some finitely generated denominator set. 
Denote $\Gamma=\{q^k: k\in\Zb\}$.
Put $F=\Fract(C)=K(q)$.\\
{\bf Definition 2.5}. Let $R$ be an unital domain,  an algebra over 
$C$ and a free $C$-module.
 Let $x$ be an element in $R$.\\
1) An element $x\in R$ is a FA-element if it is a FA-element in
$R_F:=R\otimes _CF$;\\ 
2) We say $x$ is a $\FAq$-element in $R$ if
it is a FA-element in $ R_F:=R\otimes _CF$
and for any $a\in R$ one can choose the polynomial $f_a(t)$ obeying (2.1)
such that it splites and all its roots belong to $\Gamma$.\\
{\bf Definition 2.6}. We say that two elements $a,b$ $q$-commute if 
$ab=q^kba$ for some integer $k$.\\
{\bf Proposition 2.7}[C, Prop.2.1-2.3].
Let $R$ be as in Definition 2.5 and 
$T_F=R_F[x;\tau,\delta]$ be skew extension where
$\tau$ is an automorphism, $\delta $ is a locally nilpotent 
$\tau$-derivation and $\tau\delta=q^s\delta\tau$ with $s\ne 0$. 
Denote 
$$
\wa=\sum_{n=0}^{+\infty}\frac{(1-q^s)^{-n}}{(n)_{q^s}!}
\delta^n\tau^{-n}(a)x^{-n},
\eqno(2.3)
$$
where 
$(n)_{q^s}=\frac{q^{sn}-1}{q^s-1}.$
Then \\
 1) the set $S_x=\{x^m\}_{m\in\Nb}$ is a denominator subset
 in $T_F$,\\
 2) the map $a\mapsto \wa$ is an embedding of $R$ is $T_FS_x^{-1}$,\\
 3) $T_FS_x^{-1}=\wR_F[x^{\pm 1};\tau]$ where $\wR_F$ is 
 the image of $R$
  under $a\mapsto\wa$.

  Throughout this paper $\eps$ is a primitive $l^{\rm{th}}$ root of 
1 such that $C$ admits
specialisation by $\eps: C\to K$ 
with $q\mapsto \eps$. 
 For any  
$\eps$ consider the specialisation $R_\eps$ of $R$ over $K$.
In what follows we shall use two notations.
If $a\in R$, we put  $a_\eps:=a\bmod(q-\eps)$.
For $a\in R_\eps$, we denote by $\underline{a}\in R$ an element of preimage
of $a$ under the map $\pi_\eps:R\to R_\eps=R\bmod(q-\eps)$. 
For any algebra $A$ of $R$, we denote 
$A_\eps:=(A+R(q-\eps))\bmod(q-\eps)=\pi_\eps(A)$.
  
If $u_\eps=u\bmod(q-\eps)$
lies in the center $Z_\eps$ of $R_\eps$, then
$\DD_u(a)=\frac{u\ua-\ua u}{q-\eps}\bmod(q-\eps)$ 
defines a derivation of $R_\eps$.
We call $\DD_u$ the quantum adjoint action of $u$ (see [DCKP1-2],[P3]).
An ideal is stable with respect to the quantum adjoint action 
(call $\DD$-stable) if it is
stable with respect to all $\DD_u$.
The formula $\{a,b\}=\DD_{\underline{a}}(b)$, for
$a,b\in Z_\eps$, defines the Poisson bracket on 
$\MM=\Maxspec~Z_\eps$.

Here are two versions of reduction of Proposition 2.7 modulo $q-\eps$. \\
{\bf Corollary 2.8.} Let $T,R,\tau,\delta,q^s$ be as in Proposition 2.7.
Suppose that $R$ is generated by the elements
$x_1,\ldots, x_n$ and $\tau$ is a diagonal 
automorphism with eigenvalues in $\Gamma$.
Choose $N$ such that $\delta^N(x_i)=0$ for all $1\le i\le n$.
Suppose that $l$ is relatively prime with $s$ and $l\ge N$. 
Then \\
1) $T_\eps S_{x_\eps}^{-1}\cong R_\eps[x_\eps^{\pm 1};\tau]$,\\
2) $x_\eps^l$ lies in the center $Z(T_\eps)$.\\
{\bf Proof}. 
Denote by ${\frak N}_1$ the denominator subset in $C$
generated by $q^{sn}-1$, $1\le n\le d$.
The elements $x^{\pm 1},\wx_1,\ldots,\wx_M$ 
generate $\wT:=TS_x^{-1}{\frak N}_1^{-1}$.
We denote by $\wR$ the subalgebra generated by  
$\wx_1,\ldots,\wx_M$
over $C\frak{N}_1^{-1}$.
By Proposition 2.7, the map $a\mapsto\wa$ provides
isomorphism of $R\frak{N}_1^{-1}$ onto $\wR$.
We have $\wT=\wR[x;\tau,\delta]$.  After reduction modulo $q-\eps$ we 
obtain 1).

Since $x\wx_j=q^{n_j}\wx_jx$ for some $n_j$, then $x_\eps^l$ lies in 
the center $Z(\wT_\eps)$. This proves 2).$\Box$\\
{\bf Corollary 2.9}. Let $T,R,\tau,\delta,q^s$ be as above
and $l$ be relatively prime to $s$.
Suppose that $x_\eps^l\in Z(T_\eps)$.
Then $T_\eps S_{x_\eps}^{-1}\cong R_\eps[x_\eps^{\pm 1};\tau]$.\\
{\bf Proof}.
Taking 
$$x^la=\tau^l(a)x^l + \sum^{l-1}_{i=1}{l\choose i}_{q^s}
\tau^{l-i}\delta^i(a)x^{l-i}+\delta^l(a)
$$
modulo $q-\eps$, we obtain 
$x^la = ax^l+\delta^l(a)\bmod (q-\eps)$ and
$\delta^l(a)\in(q-\eps)R$ for any $a\in R$.
If $n=lm+r$, $0\le r<l$, then
$\delta^n(a)\in(q-\eps)^mR$. On the other hand,
$(n)_{q^s}!=(q-\eps)^mc(q)$ where 
$c(\eps)\ne 0$. Hence
$$\frac{\delta^n(a)}{(n)_{q^s}!}\in Rc^{-1}(q).$$
Consider the denominator subset $\Ng_x$ in $C$
generated by $q^n-1$ where $l$ does not divide $n$ and
$\frac{q^{lm}-1}{q-\eps}$, $m\in\Nb$.
For any $a\in R$ the element $\wa$ (see 2.3) lies in the localization  of $T$ 
over $S_x$ and $\Ng_x$ and $TS_x^{-1}\Ng_x^{-1}=R\Ng_x^{-1}[x;\tau]$.
Taking modulo $q-\eps$, we get the claim. $\Box$

Let $\Sb=(s_{ij})$ be a $M\times M$ integer skew-symmetric matrix.
Denote $q_{ij}=q^{s_{ij}}$ and form the matrix $\Qb=(q_{ij})$.
Choose the subset, call distinguished subset, 
${\kg}:=\{t_1,\ldots, t_m\}$ where $1\le t_1< \ldots < t_m\le M$.\\
{\bf Definition 2.10}.
We say that $R$ is a normal quantum solvable algebra (or a NQS-algebra)
over $C$, if $R$ 
is generated by the elements $x_i$, $1\le i\le M$ and
$x_j^{-1}$, $j\in\kg$
such that the monomials 
$x_1^{t_1}\cdots x_{M}^{t_{M}}$ with $t_j\in\Zb$, $j\in \kg$ and 
$t_j\in\Nb$, $1\le j\le M$, $j\notin\kg$
form a free $C$-basis, the algebra $C$ lies in the center of $R$ and
the following
relations hold\\
1)  $x_ix_j = q_{ij}x_jx_i$ for all $i$ and $ j\in \kg$ ;\\
2) for $1\le i<j\le M$,    
$$x_ix_j=q_{ij}x_jx_i+r_{ij}\eqno (2.4)$$
 where $r_{ij}$ is a sum of monomials   
$c x_{i+1}^{t_{i+1}}\cdots x_{j-1}^{t_{j-1}}$ with 
$c\in C$. The definition of quantum solvable algebra is given in 
Remark 2.12. 

The subalgebra $Y_\kg$, generated by $C$ and $x_i^{\pm 1}$,
$i\in \kg$, is an algebra of twisted Laurent polynomials.
The subalgebras $R_i$, generated by $C$,
$x_j$, $j\ge i$ and their inversies for the distiguished 
subscripts, form a chain
$R=R_1\supset R_2\supset\cdots\supset R_M$ (call it the right filtration).
One can prove that each $R_i$ is a skew extension of $R_{i+1}$ 
[GL1,1.2].
This means that  
  the map 
$\tau_i: x_j\mapsto q_{ij}x_j, i<j$ is extended to an automorphism of
$R_{i+1}$ and the map $\delta_i:x_j\mapsto r_{ij}$ is extended to 
a $\tau_i$-derivation of $R_{i+1}$. 
All automorphisms $\tau_i$ are identical on $C$ and all 
$\tau_i$-derivations $\delta_i$ are equal to zero on $C$.
The formula (2.4) yields, $R_i=R_{i+1}[x_i;\tau_i,\delta_i]$ for 
$i\notin\kg$ and
$R_i=R_{i+1}[x_i^{\pm 1},\tau_i]$ for $i\in \kg$.
A NQS-algebra is a Noetherian domain [MC-R, 1.2.9],
a  $C$-algebra and a free $C$-module.

The NQS-algebra $R$ has the other filtration
(call it the left filtraton)
$$R'_1 \subset R_2'\subset \cdots\subset R'_{M}=R$$
with $R'_i$ is generated by $C$, $x_1,\ldots,x_i$ and their 
inversies for distinguished subscripts.
Again $R'_{i}=R_{i-1}'[x_i;\tau'_i,\delta_i']$ 
(resp. $R'_i=R'_{i-1}[x_i^{\pm 1},\tau'_i]$ for distinguished $i$)
where
$\tau_i'$ (resp. $\delta'_i$) is the automorphism 
(resp.$\tau'_i$-derivation) of
$R_{i-1}$.
We put $\delta_i=\delta_i'=0$ for distiguished $i$.

Furthermore, 
for any $1\le\alpha<\beta\le  M$ we denote by $R_{[\alpha,\beta]}$ 
the subalgebra
generated by $C$, $x_i$ and $x_j^{-1}$ such that $\alpha\le i,j\le\beta$ and 
 $j\in\kg$.

Notice that 
$R_{[\alpha,\beta]}= 
R_{[\alpha-1,\beta]}[x_\alpha;\tau_\alpha,\delta_\alpha]$,
for $\alpha\notin\kg$, and 
$R_{[\alpha,\beta]}=R_{[\alpha-1,\beta]}[x_\alpha^{\pm 1};\tau_\alpha]$,
for $\alpha \in\kg.$
Similarly,
$R_{[\alpha,\beta]}= 
R_{[\alpha,\beta-1]}[x_\beta;\tau'_\beta,\delta'_\beta]$,
for $\beta\notin\kg,$ and 
$R_{[\alpha,\beta]}=R_{[\alpha,\beta-1]}[x_\beta^{\pm 1};\tau'_\beta]$,
for $\beta \in\kg.$
We put the following conditions on a NQS-algebra.\\
{\bf Condition CN1}. We require 
$R$ be an iterated $q$-skew extension for the left and the right
filtrations.
This means that $\tau_i\delta_i=q_i\delta_i\tau_i$, for some
$q_i=q^{s_i}$, $s_i\in\Zb$,
and $\tau'_i\delta'_i=q_i'\delta'_i\tau'_i$, for some
$q_i'=q^{s_i'}$, $s'_i\in\Zb$. We require that
all $s_i\ne 0$ (resp.$s_i'\ne 0$) if $\delta_i\ne 0$ 
(resp. $\delta'_i\ne 0$). We call $\{s_i\},\{s_i'\}$
the systems of exponents $R$.\\ 
{\bf Condition CN2}. 
All $\tau_i$ and $\tau_i'$ are extended to
diagonal automorphisms of $R$ and generate 
the commuting diagonal subgroups $H$ and 
$H'$.\\
{\bf Proposition 2.11}. Let $R$ be a NQS-algebra over $C$.
Put $n=M-m$.
Let $\tR_i$, $i\notin \kg$ be a subalgebra generated by $R_i$ and $Y_\kg$. 
The chain 
$R=\tR_1\supset\tR_2\supset\cdots\supset \tR_n\supset \tR_{n+1}=Y_\kg$
is a chain of skew extensions 
$\tR_i\cong\tR_{i+1}[x;\ttau_i,\tdelta_i]$.
If, in addition, $R$ obeys Condition CN1, then
$\ttau_i\tdelta_i=q_i\tdelta_i\ttau_i$
with the same $q_i=q^{s_i}$ as in CN1.\\
{\bf Proof}. We put $\ttau_i(a)=\tau_i(a)$ (resp.$\tdelta(a)=\delta(a)$), 
for $a\in R_{i-1}$,
and $\tau_i(x_j)=q_{ij}x_j$ (resp.$\tdelta(x_j)=0$), for $j<i$, $j\in\kg$.
The direct calculations conclude the proof.$\Box$\\
{\bf Remark 2.12}. A quantum solvable algebra is defined in [P1-P3] 
as an algebra generated  
$R$ is generated by the elements 
$x_1, x_2,\ldots, x_n, x_{n+1}^{\pm 1},\ldots,
x_{n+m}^{\pm 1}$ with $M=n+m$
such that the monomials $x_1^{t_1}\cdots x_n^{t_n}x_{n+1}^{t_{n+1}}\cdots
x_{n+m}^{t_{n+m}}$ with $t_1,\ldots,t_n\in\Nb$ and 
$t_{n+1},\ldots,t_{n+m}\in\Zb$ 
form a free $C$-basis and the relations hold:
1) $x_ix_j=q_{ij}x_jx_i$, for all $i$ and $n+1\le j\le M$,\\
 2) 
$x_ix_j=q_{ij}x_jx_i+r_{ij}$, $1\le i<j\le n$
where $r_{ij}$ is an element of the subalgebra
$R_{i+1}$ generated by $x_{i+1},\ldots, x_n, x_{n+1}^{\pm 1},\ldots,
x_{n+m}^{\pm 1}.$.
Proposition 2.9 claims that a NQS-algebra is a quantum solvable algebra.
The Conditions CN1 and CN2 are comparable with more general Conditions
Q1-Q4 of [P2] and Conditions 3.2-3.4 of [P3].\\
{\bf Proposition 2.13}. Any FA-element in a quantum solvable algebra 
(in particular, is a NQS-algebra) $R$ is
a $\FAq$-element.\\
{\bf Proof}. 
Let $R$ be a quantum solvable algebra (see above Remark).
For a monomial $w=x_1^{t_1}\cdots x_{M}^{t_{M}}$, denote
$\deg(w)=(t_1,\ldots,t_{M})$.
Lexicographical order provides the filtration in $R$.
The algebra $A_\Qb:=\gr(R)$ is 
generated by $a_i=\gr(x_i)$, $1\le i\le M$ and $a_j^{-1}$, $j\in \kg$.
The relatons are $a_ia_j=q_{ij}a_ja_i$.
The algebra $A_\Qb$ is the localization of 
algebra of twisted polynomials.
As usual $(\cdot,\cdot)$ denotes the standard scalar
product in $\Zb^M$. For two monomials 
$a,b\in A_\Qb$ with $\deg(a)=\um$, $\deg(b)=\un$,
we have $ab=q^{(\Sb \um,\un)}ba$.
For every $u,v\in R$ with $\deg(u)=\um$, $\deg(v)=\un$,
$$
uv=q^{(\Sb \um,\un)}vu+\{\mbox{terms of lower degree}\}.\eqno(2.5)
$$
Let $u,v\in R$ and the element $u$ be a FA-element.
Let $f(t)$ be the corresponding polynolial obeying (2.1), for 
$x:=u$ and $a:=v$.
Put $\gamma:=q^{(\Sb\um,\un)}.$
By (2.5),
$$0 = c_0u^Nv + u^{N-1}vu + \ldots + c_Nvu^N =
f(\gamma)vu^N+\{\mbox{terms of lower degree}\}.$$
Hence, $f(\gamma)=0$ and $f(t)=(t-\gamma)f_1(t)$.
The element $v_1:=uv-\gamma vu$ is annihilated by
$f_1(\Ad_u)$. The proof is concluded by induction on degree
of polynomial $f(t)$.$\Box$

Here are two the most familiar examples of NQS-algebras.\\   
{\bf Example 2.14}. Quantum matrices.\\
The algebra $M_{q}(n,K)$ of regular functions on quantum matrices
is generated over $C:=K[q,q^{-1}]$ and the entries of quantum matrix  
$\{a_{ti}\}_{t,i=1}^n$ which obey the relations 
$a_{ti}a_{sj}- a_{sj}a_{ti}= 
(q-q^{-1})a_{si}a_{tj}$
for $i<j$, $t<s$  and
$ a_{ti}a_{sj}= qa_{sj}a_{ti} $, for $t<s, i=j$ and $t=s, i<j$. 

The  algebra $M_q(n,K)$ is a NQS-algebra
with respect to generators $x_{(i-1)n+j}=a_{ij}$.
It obeys CN1 (see [G],[P2]) and CN2 
(the map $\tau_{ij}:a_{ij}\mapsto qa_{ij}$,
i.e the multiplication of $i$th row by $q$, is an
automorphism of $R$).  
One can obtain an other examples
considering subalgebras (like Quantum triangular matrices),
some generalizations and muliparameter versions of this algebra.\\
{\bf Example 2.15}. $U_q(\nog)$, where $\nog$ is the upper nilpotent subalgebra
of semisimple Lie algebra. The algebra $U_q(\nog)$
is generated over $ C=K[q,q^{-1}, (q^{d_i}-q^{-d_i})^{-1}]$
by $E_i$, $i=\overline{1,n}$ with the quantum 
Chevalley-Serre relations.
Fix a reduced expression $w_0=s_{i_1}\ldots s_{i_N}$ of the longest element
in the Weyl group W. Consider the following convex ordering 
$\beta_1=\alpha_{i_1}, \beta_2=s_{i_1}(\alpha_2),\ldots, 
\beta_N= s_{i_1}\ldots s_{i_{N-1}}(\alpha_N)$
in the set $ \Delta^+$ of positive roots.
Consider the quantum root vectors 
$ E_{\beta_s}=T_{i_1}\cdots T_{i_{s-1}}E_{i_s}$, for $1\le s\le N$ 
  [Lu].
There are following relations on the $E_{\beta_i}$'s [LS]:
$$E_{\beta_i} E_{\beta_j}- q^{-(\beta_i,\beta_j)}
E_{\beta_j}E_{\beta_i}= \sum_{m\in\Zb_+^N}c_m E^m,$$
where $i<j$, $c_m\in K[q,q^{-1}]$ and $c_m\ne 0$ only 
when $m=(m_1,\ldots,m_N)$ is such that $m_s=0$ for $s\le i$ and $s\ge j$.
The algebras $U_q(\nog)$, $U_q(\bg)$, and subalgebras $U^w_q(\nog)$ 
(see [C1],[DC-P1]) are NQS-algebras.
They obey Conditions CN1 (see [G], [P2]) and CN2
(the map $\tau_\alpha:E_\beta\mapsto q^{(\alpha,\beta)}E_\beta$ is an
automorphism of $R$).
\\
{\bf Proposition 2.16.} Let $R$  be a NQS-algebra obeying Condition CN1 with
the systems of exponents $\{s_i\},\{s'_i\}$. \\
1) All $x_\alpha$, $1\le \alpha\le M$ are FA-elements in $R$.
Choose $N_\alpha$ (see Proposition 2.4) such that
$\delta_\alpha^{N_\alpha}(x_j)=0$, 
$\alpha<j$ and $(\delta'_\alpha)^{N_\alpha}(x_j)=0$, $j<\alpha$.\\ 
2) For any  $1\le \alpha\le M$, $\alpha\notin \kg$
consider two denominator subsets
$\Ng_\alpha$ generated by $q^{s_\alpha t}-1$,
$q^{s'_\alpha t}-1$, $1\le t<N_\alpha$ and 
$S_\alpha$ generated by $x_\alpha$. The algebra 
$RS_\alpha^{-1}\Ng_\alpha^{-1}$ is a NQS-algebra with distinguished subset 
$\kg\bigcup \{\alpha\}$ over $C\Ng_\alpha^{-1}$, with the same (as $R$)
matrix $\Qb$ and 
systems of exponents.\\
{\bf Proof}.
The claim 1) is proved similarly [P1,Lemma 4.3].
To prove 2) we apply Proposition 2.5 for two extensions
$R_\alpha=R_{\alpha+1}[x_\alpha;\tau_\alpha,\delta_\alpha]$
and $R'_\alpha=R'_{\alpha-1}[x_\alpha;\tau'_\alpha,\delta'_\alpha]$, 
and
consider the new system of generators
$\wx_1,\ldots,\wx_{\alpha-1},x_\alpha^{\pm 1},\wx_{\alpha+1},\ldots,\wx_M$
of $RS_\alpha^{-1}$.$\Box$.\\
{\bf Corollary 2.17}. 
Let $R$ be as in Proposition 2.16. Suppose that $l$ is relatively prime with
$s_\alpha$,$s'_\alpha$
and  $x_{\alpha\eps}^l$ lies in the center of $R_\eps$,
then $R_\eps S_{\alpha\eps}^{-1}$ is isomorphic to the specialisation of some  
NQS-algebra modulo $q-\eps$.\\
{\bf Proof}.
Consider the multiplicatively closed subset
$\Ng_{\alpha,l}$ generated by $q^{s_\alpha n}-1$, $q^{s'_\alpha n}-1$
(where $1\le n < N_\alpha$ and $l$ does not divide $n$) 
and $\frac{q^{s_\alpha lm}-1}{q-\eps}$, $\frac{q^{s'_\alpha lm}-1}{q-\eps}$
for 
$1\le lm < N_\alpha$.
Since $l$ is relatively prime with $s_\alpha$ and $s'_\alpha$, 
then polynomials of $\Ng_{\alpha,l}$
are not zero at $q=\eps$. The element $x_{\alpha\eps}^l$ lies in the center
of $R_\eps$; by the proof of Corollary 2.9,   
$\wx_1,\ldots,\wx_{\alpha-1},x_\alpha^{\pm 1},\wx_{\alpha+1},\ldots,\wx_M$ lie
in $RS_\alpha^{-1}\Ng_{\alpha,l}^{-1}$. One can reduce the generators modulo
$q-\eps$ and get the system
of generators $\wx_{1\eps},\ldots,\wx_{\alpha-1,\eps},x_{\alpha\eps}^{\pm 1},
\wx_{\alpha+1,\eps},\ldots,\wx_{M\eps}$ of
$R_\eps S_{\alpha\eps}^{-1}$.

For a NQS-algebra $R$, consider $N:=N_R=\rm{max}\{N_\alpha\}$.
For $1\le i_1<\ldots<i_k\le M$, $\mu:=\{i_1,\ldots,i_k\}\supset\kg$
we denote by $\Sb_\mu$ the submatrix $(s_{ij})$, $i,j\in\mu$ of $\Sb$.\\
{\bf Definition 2.18}. 
We say that a positive integer $l$ (resp. a primitive 
$l^{\mbox{th}}$ root of unity
$\eps$)
is admissible for a NQS-algebra $R$ 
if it obeys the conditions:\\
1) $l$ is relatively prime with all elementary divisors
of all submatrices $\Sb_\mu$, $\mu\supset\kg$;\\
2) $l$ is relatively prime with  $s_i,s_i'$, $1\le i \le M$;\\
3) $l\ge N$.\\
{\bf Lemma 2.19}. Let $\eps$ obeys the conditions 2) and 3)
of Definition 2.18, and $R$ be a NQS-algebra 
 obeying Condition CN1.
Then
the elements  $\{x_{i\eps}^l\}$ lie in the center $Z_\eps$ of $R_\eps$.\\
{\bf Proof}. Apply Corollary 2.8.$\Box$\\
{\bf Proposition 2.20}. Let $R$ and $\eps$ be as in 2.19.\\
1) If $x$ is a FA-element of $R_F$ (resp. $R_\eps$), 
then linear operator $\Ad_x$ is diagonalizable in $R_FS_x^{-1}$ (resp. 
$R_\eps S_x^{-1}$).\\
2) For any FA-element $x$ in $R$ the element $x_\eps^l$ lies in $Z_\eps$. \\
{\bf Proof}.
Lemma 2.19 implies that $R_\eps$ is finite over its center. 
The set of roots of unity, 
that obey 2) and 3) of Definition 2.18, is infinite.
The statement 1) is a corollary  [P1, Cor.2.5, Proposition 3.4].

 Let us prove 2). 
 For the FA-element $x$ and any 
$a$ is $R$ there exists a minimal polynomial $f(t)$ that obeys (2.1).
By Proposition 2.13, the roots of $f(t)$ belong to $\Gamma$.
Suppose that $q^{\alpha_1},\ldots,q^{\alpha_N}$ are the roots of $f(t)$.
The element $u=x^l$ are also a FA-element of $R$.
The operators $\Ad_x$ and $\Ad_{x^l}$ are simultaneously diagonalizable.
The roots of corresponding polynomial $f_*(t)$ for $x^l$ are
$\lambda_i:=q^{\alpha_il}$, $1\le i \le N$. 
It implies $f_*(t)=(t-1)^N\bmod(q-\eps)$
and $(\Ad_{x_\eps^l}-\mbox{id})a=0$.
On the other hand, by 1), $x_\eps$ and $x_{\eps}^l$ are FA-elements
in $R_\eps$. Hence, $\Ad_{x_\eps^l}$ is diagonalizable.
It follows $x_\eps^l\in Z_\eps$ and 2).$\Box$\\
{\bf Definition 2.21}. $R$ and $\eps$ as above.
We say that an ideal is $\DD_0$-stable if it is stable with respect to
all derivations $\DD_{x_i^l}$, $1\le i\le M$.\\
{\bf Notation 2.22}. For any automorphism $\tau\in H$ we denote by
$\theta$ the following diagonal derivation of $R_\eps$
$$\theta(a)=\frac{\tau^l-{\rm {id}}}{q-\eps}~\bmod(q-\eps).$$
Similar $\theta'$ for $\tau'\in H'$.
By $\Theta$ we denote the commutative subalgebra spanned
by $\theta_1,\ldots,\theta_M$. Similarly for $\Theta'$.

\section{Stratification of prime ideals}
In this section, we stratify the 
prime spectrum of $R$ and the prime $\DD$-stable spectrum of $R_\eps$
(Theorem 3.2). It is proved that any prime $\DD$-stable ideal of $R_\eps$
is completely prime (Theorem 3.3). 

Throughout this section $R$ is a NQS-algebra,
obeying Conditions CN1 and CN2,
and $\eps $ obeys the conditions 2) and 3)
of Definition 2.18.

Consider the mulplicatively closed subset
$\Ng=\prod_\alpha \Ng_{\alpha,l}$ (see Corollary 2.17).
The polynomials of $\Ng$ are not zero at $q=\eps$.

Fix an integer $i_1$ which $1\le i_1\le M$.
If $ i_1\in \kg$, we put $R^{(1)}:=R$. If $i_1\notin\kg$, 
we consider the denominator subset 
$S_1$ generated by $y_1:=x_{i_1}$.
According to Proposition 2.16 and Corollary 2.17,  
$R^{(1)}:=RS_1^{-1}\Ng^{-1}$ is a NQS-algebra
over $C\Ng^{-1}$ with the same (as $R$) systems of exponents.
The algebra $R^{(1)}$ is generated by
$$
x'_1,\ldots, x'_{i_1-1},x_{i_1}^{\pm 1},x'_{i_1+1},
\ldots,x'_{M}\eqno (3.1)
$$
where 
$x'_j:=\wx_j$.
Recall that all generators
$q$-commute with $y_1$ and 
are FA-elements in $R^{(1)}$.
It follows that, for all $i$, the elements ${x'_{i\eps}}^l$ lie 
in the center of 
$R^{(1)}_\eps= R_\eps S_{1\eps}^{-1}$ (see Proposition 2.20).

Let $i_2$ be any integer which $i_1< i_2\le n$.
There exists a positive integer $t$ such that 
$y_2:=x'_{i_2}x_{i_1}^t\in R$.
Similarly to the first step of stratification process,
we consider denominator subset $S_2$ generated by 
$q$-commuting elements $y_1,y_2$. 
As we saw the element
$(x'_{i_2\eps})^l$ lies in the center of $R^{(1)}_\eps$.
By Corollary 2.17, 
the algebra $R^{(2)}:=RS_2^{-1}\Ng^{-1}= RS_1^{-1}S_{x_{i_2}'}^{-1}\Ng^{-1}$
 is a NQS-algebra 
with the generators
$$
x''_1,\ldots, x''_{i_1-1},x_{i_1}^{\pm 1},x''_{i_1+1},
\ldots, x''_{i_2-1}, {x'_{i_2}}^{\pm 1},x''_{i_2+1}\ldots, x''_M.\eqno (3.2)
$$

After $k$ steps we get
the denominator subset $S:=S_\mu$, $\mu:=\{i_1,\ldots,i_k\}$ 
 generated by the  
system of $q$-commuting elements $y_1,\ldots,y_k\in R$ and 
$\Ng$. We call $S$ as the standard denominator subset.
The algebra $\tR:= R^{(k)}= RS^{-1}$ is a NQS-algebra over $C\Ng^{-1}$
with the generators $\tx_j:=x_j^{(k)}$ and $y_1^{\pm 1},\ldots, y_k^{\pm 1}$.
All generators are FA-elements in $\tR$.

We denote by $Y:=Y_\mu$ the subalgebra, 
generated by $y_1^{\pm 1},\ldots, y_k^{\pm 1}$
(or $x_{i_1}^{\pm 1}$, $(x'_{i_2})^{\pm 1}$, $(x''_{i_3})^{\pm 1}$,
$\ldots, (x^{(k-1)}_{i_k})^{\pm 1}$).
The imposed relations of $Y$ are $y_iy_j=q^{t_{ij}}y_jy_i$.
The integer matrix $(t_{ij})_{i,j=1}^k$ is obtained by elementary
transformations of submatrix $\Sb_\mu=(s_{ij})_{ij\in\mu}$ of $\Sb$.
The algebra $Y$ is an algebra of twisted Laurent polynomials.

By Proposition 2.11, we may treat $\tR$ as a iterated $q$-skew extension 
$$\tR:=
\tR_1\supset\tR_2\supset\cdots\supset
\tR_\tM\supset\tR_{\tM+1}=:Y\eqno(3.3)$$
where $\tM:=M-k$ and $\tR_i\cong\tR_{i+1}[x;\ttau_i,\tdelta_i]$.\\
{\bf Definition 3.1.}\\
1) We say that $S:=S_\mu$ is $C$-admissible  
if the ideal $\JJ:=\JJ_S$ of $\tR$ generated by 
$\tx_i$, $i\in[1,M]-\mu$ has zero intersection with $C$;\\
2) We say that $S:=S_\mu$ is $\eps$-admissible 
if the ideal $J:=J_S$ of $\tR_\eps$  generated by 
$\tx_{i\eps}$, $i\in[1,M]-\mu$ is proper.\\
3) We say that $S:=S_\mu$ is $(\eps,\DD)$-admissible 
if the $\DD$-stable ideal (denote $\DDJ$ ) 
of $\tR_\eps$ generated by $J$ (see 2))
is proper. 

Notice that, in general, the ideal $\JJ$ (resp. $J$) 
may have nonzero intersection with
$Y$ (resp. $Y_\eps$) and is not prime .
For instance, it holds for the algebra $R_f$ which is constructed by a
 polynomial $f$ as follows. This algebra is
generated by $x_1, x_2, y^{\pm 1}_1,\ldots,
y_k^{\pm 1}$ where the elements $\{y_i\}$ lie in the center and
$x_1x_2-qx_2x_1=f(y_1,\ldots, y_k, q)$. The ideal $\JJ$, generated by $x_1$
and $x_2$, has nonzero intersection with $Y$ and is not prime when
the polynomial $f$ is reducible. In the case $f=f(q)$, 
the ideal $\JJ$ has nonzero  intersection with $C$.\\
{\bf Theorem 3.2}. Let $R$ be a NQS-algebra obeying  
Conditions CN1 and CN2, 
and 
$\eps$ be a specialisation of $C$ obeying the conditions 2) and 3)
of Definition 2.18. \\
1) For any $\II\in\Spec(R)$, $\II\bigcap C=0$, there
exists a unique $C$-admissible standard denominator subset 
$S:=S_\mu$ such that   
$\II\bigcap S=\emptyset$ and $\II S^{-1}\supset\JJ_S$.\\
2) Let $R$ and $\eps$ be as above.
For any prime $\DD$-stable ideal $I$ of $R_\eps$ there
exists a unique $(\eps,\DD)$-admissible  standard denominator subset 
$S=S_\mu$
such that 
$I\bigcap S_\eps=\emptyset$ and $IS_\eps^{-1}\supset J_S$.\\
{\bf Proof}. 
Let $\II\in\Spec(R)$ and $\II\bigcap C=0$.
Suppose that $x_1\ldots,x_{i_1-1}\in \II$ and $y_1:=x_{i_1}\notin \II$.
All prime ideals of $R$ are completely prime [GL2,Theorem 2.3].
(this is false for $R_\eps$).
Therefore, $\II\bigcap\{y_1^t\}_{t\in\Nb} = \emptyset$.
The ideal $\II$ admits localization over $S_1\cdot \Ng$ 
(see stratification process).
By the formula (2.3), $x_1',\ldots,x'_{i_1-1}\in\II S_1^{-1}\Ng^{-1}$.
Suppose that $x'_ {i_1+1},\ldots,x'_{i_2-1}\in\II S_1^{-1}\Ng^{-1}$ and 
$x_{i_2}'\notin \II S_1^{-1}\Ng^{-1}$.
Following the stratification process, finally, we have $\tx_i\in\II S^{-1}$ for
$i\in[1,M]-\mu$. This proves 1).

 For a prime $\DD$-stable ideal $I$ of $R_\eps$, consider the greatest 
$\Theta'$-stable ideal $I_{\Theta'}$ in $I$. The ideal $I_{\Theta'}$ 
is $\DD$-stable
[P3, Proposition 3.14] and prime [MC-R,14.2.3],[D,3.3.2].
Consider the left filtration 
$R'_{1\eps}\subset\cdots\subset 
R'_{i\eps}\subset\cdots\subset R'_{M\eps}=R_\eps$.
A prime ($\DD$, $\Theta'$)-stable ideal has prime intersections with all
subalgebras $R'_{i\eps}$ [P3, Theorem 2.12].
Suppose that $I$ contains
$x_{1\eps},\cdots,x_{i_1-1,\eps}$ and does not contain
 $y_{1\eps}=x_{i_1,\eps}$.
The ideal $I_{\Theta'}\bigcap R'_{i\eps}$ is prime. 
Since
$$
\frac{R'_{i\eps}}{I_{\Theta'}\bigcap R'_{i\eps}}\cong 
\frac{K[x_{i\eps}]}{I_{\Theta'}\bigcap K[x_{i\eps}]},$$
the ideal $I_{\Theta'}\bigcap R'_{i\eps}$ is completely prime.
It follows that 
$I_{\Theta'}$ has 
empty intersection with the subset 
$S_{1\eps}:=\{y_{1\eps}^m\}_{m\in \Nb}$.
Since $y_{1\eps}$ is a $\Theta'$-eigenvector,
$I$ has empty intersection with $S_{1\eps}$.

Since $y_1$ is a FA-element in $R$, the element
$y_{1\eps}^l$ lies in the center of $R_\eps$ (see Proposition 2.20).
By proof of Corollary 2.9, we can reduce $x'_1,\ldots,x'_{i_1-1}$
modulo $q-\eps$ and obtain $x'_{1\eps},\ldots,x'_{i_1-1,\eps}\in 
IS_{1\eps}^{-1}$.

On the second step suppose that
$x'_ {i_1+1,\eps},\ldots,x'_{i_2-1,\eps}\in
I_{\Theta'}S_{1\eps}^{-1}$ and $x'_{i_2}\notin I_{\Theta'}S_{1\eps}^{-1}$.
The factor algebra  
$$\frac{R'_{i_2}S_{1\eps}^{-1}}{I_{\Theta'}\bigcap R'_{i_2}S_{1\eps}^{-1}},
$$ is a prime factor of
the algebra of generated by two $q$-commuting $y_1$ and $x'_{i_2}$.
The image of $x'_{i_2}$ is either zero or regular [P3,Lemma 3.11].
Since $x'_{i_2}\notin I_{\Theta'}$, the image is regular.
The ideal $I_{\Theta'}$ (and $IS_{1\eps}^{-1}$) has empty intersection
with $S_{2\eps}$ generated by $y_{1\eps}$ and $y_{2\eps}$ 
(see stratification process).
We consider localization over
$S_{2\eps}$. After $k$ steps we obtain 2).
$\Box$

We say that an ideal of $Y_\eps$ is $\DD_0$-stable,
if it is  stable with respect to all derivations
$\DD_{y_i^l}$, $1\le i\le k$.
By direct calculations [P3,Lemma 3.16],
$$\DD_{y_i^l}y_{j\eps}=t_{ij}l\eps^{-1}y_{j\eps}y_{i\eps}^l,$$
$$\{y_{i\eps}^l,y_{j\eps}^l\}=t_{ij}l^2\eps^{-1}y_{i\eps}^ly_{j\eps}^l.$$
{\bf Theorem 3.3}. Let $R$, $\eps$ be as Theorem 3.2.
Any prime $\DD$-stable ideal of $R_\eps$ is completely prime.\\
{\bf Proof}. 
Let $I$ be a prime $\DD$-stable ideal of $R_\eps$.
According the previous Theorem,
$$
\frac{R_\eps S_\eps^{-1}}{IS_\eps^{-1}}=
\frac{Y_\eps}{IS_\eps^{-1}\bigcap Y_\eps}\eqno(3.4)
$$
where $Y_\eps$ is the algebra of twisted Laurent
polynomials generated by $y_{1\eps},\ldots,y_{k\eps}$.
It follows that the ideal $IS_\eps^{-1}\bigcap Y_\eps$ of $Y_\eps$ is prime.
 Since $y_i$ is a FA-element, $y_{i\eps}^l$ lies
in the center $Z_\eps$ (see Proposition 2.20).
The ideal $I$ is $\DD$-stable, hence, it is stable with respect to 
$\DD_{y_i^l}:R_\eps\to R_\eps$, $1\le i\le k$.
The same is true for  $IS_\eps^{-1}\bigcap Y_\eps$.
Any prime $\DD_0$-stable
ideal of an algebra of twisted Laurent polynomials
is completely prime [P3, Corollary 3.18].
The ideal  $IS_\eps^{-1}\bigcap Y_\eps$ is completely prime  and, 
therefore, $I$ is completely prime. $\Box$

Till the end of this section we suppose 
that $\eps$ is an admissible specialisation of $C$
(see Definition 2.18).
 One can choose the new generators (monomials)  
$h_1,\cdots,h_k$ of $Y$ with the following relations
$$h_1h_2=q^{m_1}h_2h_1, \ldots,h_{2r-1}h_{2r}=q^{m_r}h_{2r}h_{2r-1}\eqno (3.5)
$$
where $m_1,\ldots,m_r$ are relatively prime with $l$ 
(see Definition 2.18)
and $h_{2t+1},\ldots,h_k$ generate the center of $Y$.
All elements $h_i$ are FA-elements in $\tR$.
In what follows we suppose that the elements of $C$ and 
$z_1:=h_{2r+t+1},\ldots,z_p:=h_k$, $p=k-2r-t$
generate the intersection $\Zg:=Y\bigcap\tZ$ where $\tZ:=\rm{Center}(\tR)$.
Denote $u_1:=h_{2r+1},\ldots, u_t:=h_{2r+t}$.
We have $\Zg=K[z_1^{\pm 1},\ldots, z_p^{\pm 1},q^{\pm 1}]$,
$\Zg_\eps=K[z_{1\eps}^{\pm 1},\ldots, z_{p\eps}^{\pm 1}]$ and 
$$ 
Z(Y)= K[u_1^{\pm 1},\ldots,u_t^{\pm 1},
z_1^{\pm 1},\ldots, z_p^{\pm 1},q^{\pm 1}],
$$
$$
Z(Y)_\eps:=Z(Y)\bmod(q-\eps)= K[u_{1\eps}^{\pm 1},\ldots,u_{t\eps}^{\pm 1},
z_{1\eps}^{\pm 1},\ldots, z_{p\eps}^{\pm 1}   ],
$$
$$
Z(Y_\eps):=\rm{Center}(Y_\eps)= K[h_{1\eps}^{\pm l},\ldots,h_{2r\eps}^{\pm l},
u_{1\eps}^{\pm 1},\ldots,u_{t\eps}^{\pm 1},z_{1\eps}^{\pm 1},
\ldots, z_{p\eps}^{\pm 1}].
$$
The algebra $Z(Y)_\eps$ coincides with subalgebra $Z(Y_\eps)^\DD$
which consists of the elements of $Z(Y_\eps)$ annihilated by all 
$\DD_{y_i^l}$. 
As above $\tZ_\eps:=\rm{Center}(\tR_\eps)$.
The intersection $\tZ_\eps\bigcap Y_\eps$ is a polynomial algebra
$$
\tZ_\eps\bigcap Y_\eps=K[h_{1\eps}^{\pm l},\ldots,h_{2r\eps}^{\pm l},
u_{1\eps}^{\pm l},\ldots,u_{t\eps}^{\pm l},z_{1\eps}^{\pm 1},
\ldots, z_{p\eps}^{\pm 1}].$$
{\bf Notations 3.4}.\\
1) $G$ is the subgroup in $\tR$ generated by $S$ (i.e. by $y_1,\ldots,y_k$),\\
2) $G^l$ is its subgroup generated by $y^l,\ldots, y_k^l$,\\
3) $W:=\{ a\in \tR: ay=ya~\mbox{ for all}~ y\in Y\}$.\\
4) $W_\eps:=W\bmod(q-\eps)$.

The elements of $G$ are FA-elements on $\tR$.
It follows that, for any $y\in G$,  the linear operator 
$\Ad_y$ is diagonalizable over $C\Ng^{-1}$, 
 (Proposition 2.20).
 Since the generators of $G$ are $q$-commuting elements,
 then $\{\Ad_y: y\in G\}$ is the commutative subgroup of $\rm{Aut}(\tR)$.
It follows that $\{\Ad_y\}$ are simultaneously diagonalizable.  
  
The map
$\Delta_{y^l}:= y^{-l}_{\eps}\DD_{y^l}:\tR_\eps\to \tR_\eps$ is a
diagonalizable derivation. Moreover, if $\Ad_y v=q^{\alpha} v$, then
$\Delta_{y^l}(v_\eps) = \ual v_\eps$ where $\ual:=\alpha l\eps^{-1}$.
If $\DD_{y^l}(v_\eps)=0$ for any $y\in G$, then $v_\eps\in W_\eps$.
\\
The derivation $\Delta_{y^l}$ 
preserve the center $\tZ_\eps$ and diagonalizable in it. 
\\
{\bf Lemma 3.5}.
  Let $v\in \tR$ and $v_\eps\in\tZ_\eps$. 
 Then\\ 
 1) $\DD_v(Y_\eps)$ is contained in the 
ideal $<v_\eps>$ generated by $v_\eps$.\\
2) If $v_\eps\in W_\eps\bigcap \tZ_\eps$, then $\DD_v(Y_\eps)=0$. \\
{\bf Proof}.
One can assume, that $v_\eps$ (resp. $v$) is $\Delta_{G^l}$-eigenvector
(resp. $\Ad_G$-eigenvector).
For $\Ad_yv=q^\alpha v$, we have $\Delta_{y^l}v_\eps=\ual v_\eps$ and
$\ual y_\eps^l v_\eps= \DD_{y_l}v_\eps$.
On the other hand,
$$\DD_{y^l}v_\eps= \{y_\eps^l, v_\eps\}=
-\{v_\eps,y_\eps^l\}=-\DD_v(y^l_\eps)=
-ly_\eps^{l-1}\DD_v(y_\eps).$$ 
It follows 
$$\DD_v(y_\eps)=-\ual l^{-1}y_\eps v_\eps.\eqno (3.6)$$
Formula (3.6) yields 1).

To prove claim 2), we decompose $v$ into the sum $v=v_0+v_1+\ldots+v_n$
of 
$\Ad_G$ eigenvectors. Suppose that $v_0\in W$ (i.e $\Ad_yv=v$ for all $y\in G $)
 and $\Ad v_i=q^{\alpha_i}v_i$, $\alpha_i\ne 0$ for $1\le i\le n$.
 Since $v_\eps\in W_\eps$, then $v_{i\eps}=0$ for $1\le i\le n$.
 Using (3.6), we have 
 $$\DD_v(y_\eps)=\DD_{v_0}(y_\eps)+\DD_{v_1}(y_\eps)+\ldots+\DD_{v_n}(y_\eps)=
 \DD_{v_0}(y_\eps)=0.$$
 This proves 2). $\Box$\\
{\bf Proposition 3.6}. Let $S:=S_\mu$ be $(\eps,\DD)$-admissible
 and $\DDJ$ denotes the lowest $\DD$-stable ideal which containes $J:=J_S$
 (see Definition 3.1).
Then
$$\tZ_\eps= \DDJ\bigcap \tZ_\eps + \Hle (W_\eps\bigcap\tZ_\eps).\eqno (3.7)$$
{\bf Proof}. Let $v_\eps\in \tZ_\eps$ be 
a common eigenvector for $\Delta_{G^l}$ .
If $v_\eps\notin\DDJ\bigcap \tZ_\eps$, then $v_\eps=j_{0\eps}+y_{0\eps}$
where $j_{0\eps}\in \DDJ$, $y_{0\eps}$ is a nonzero element
of $Z(Y_\eps)$ and $j_{0\eps}$, $y_{0\eps}$ are $\Delta_{G^l}$-eigenvectors
with the common eigenvalue.
One can present $y_{0\eps}$ in the form 
$
y_{0\eps}= h_\eps y_{0\eps}'
$
 where $h_\eps$ is some monomial :
$$h_\eps:=h_{1\eps}^{m_1l}\cdots
h_{2t\eps}^{m_{2t}l}$$
with
$m_1\ldots,m_{2t}\in\Zb$, $y_{0\eps}'\in Y_\eps$ 
and $\Delta_{G^l}y_{0\eps}'=0$.
Then the element $v_{\eps}':= h_\eps^{-1}v_\eps$ obeys 
$\Delta_{G^l}v_{\eps}'=0$. 
Whence $v_{\eps}'\in W_\eps$.$\Box$\\
 Notice that
 $\DDJ\bigcap Y_\eps$ is a $\DD_Y$-stable ideal in the algebra of 
 twisted Laurent polynomials $Y_\eps$.
 Hence [P3,Lemma 3.17], $\DDJ\bigcap Y_\eps$ is
  is generated by its intersection with $Z(Y)_\eps$.\\
{\bf Proposition 3.7}
1) If ${\frak m}$ is a maximal ideal of $Z(Y)_\eps$ which
lies over $\DDJ\bigcap Z(Y)_\eps$, then 
$L_{\frak m}:=\DDJ + {\frak m}Y_\eps$ is $\DD$-stable ideal in $\tR_\eps$;\\
2) if ${\frak M}$ is a maximal ideal of $W_\eps \bigcap \tZ_\eps$
over $\DDJ\bigcap W_\eps\bigcap \tZ_\eps$,
then
$$L_{\frak M}:=\DDJ\bigcap \tZ_\eps + \Hle {\frak M}$$
 is a Poisson ideal of $\tZ_\eps$.\\
{\bf Proof}. 
 By the formula $\tR_\eps=\DDJ+Y_\eps$ (resp. (3.7)),
 $L_{{\frak m}}$ (resp. $L_{{\frak M}}$) is a two-sided ideal in 
 $\tR_\eps$ (resp. $\tZ_\eps$).
 
Let $v_\eps\in \tZ_\eps$ (resp. $v\in \tR$)  
be a common $\Delta_{G^l}$-eigenvector
(resp. $\Ad_{G}$). We are going to prove that both ideals $L_{\frak m}$ and
$L_{\frak M}$ are $\DD_v$-stable.

If $v_\eps\in \DDJ\bigcap\tZ_\eps$ or
$v_\eps\in W_\eps\bigcap\tZ_\eps$ the statement is a corollary of Lemma 3.5.
For $v_\eps\in \Hle $, the derivation  
$\DD_v$ is zero in $Z(Y)_\eps$ and $W_\eps$.
Both ${\frak m}$ and ${\frak M}$ are annihilated by $\DD_v$.
The ideals $L_{\frak m}$ and
$L_{\frak M}$ are $\DD_v$-stable.
  $\Box$

\section{Irreducible representations}

Let $R$ be an algebra and a free $C$-module.
One can consider specialisation $R_\eps$ of $R$.
As above $Z_\eps$ is a center of $R_\eps$.
This algebra has a Poisson structure via quantum adjoint action
(see Section 2). 
 
Let $\chi$ be a central character $\chi:Z_\eps\to K$
and $m(\chi)$ the corresponding maximal ideal. We treat $\chi$ as a point of
variety $\MM:=\Maxspec(Z_\eps)$. 
We consider stratification of $\MM$ [BG1]:
$\MM=\MM_0\supset\MM_1\supset\cdots\supset \MM_m=\emptyset$
where
$\MM_{i+1}=(\MM_i)_{\rm{sing}}$.
All $\MM_i$ are Poisson varieties.
In the case $K=\Cb$, the smooth locuses 
$\MM_i^0:=\MM_i-\MM_{i+1}$ are complex analytic Poisson varieties.
Each symplectic leaf is  a disjoint union of symplectic leaves.
 For $\chi$ we denote by $\Omega_\chi$ 
the corresponding
symplectic leaf. 

Let $m(\chi,\DD)$ be the greatest Poisson (i.e. $\DD$-stable) 
ideal in $m(\chi)$.
One can treat the algebra
$\FF:=Z_\eps/ m(\chi,\DD)$ as the algebra of regular functions
on Zariski closure $\MM_\chi$ of $\Omega_\chi$.
Denote by $R_{\eps,\chi}$ the finite dimensional subalgebra
$R_\eps/m(\chi)R_\eps$.
\\
{\bf Lemma 4.1}[P3, Lemma 5.1].
Let $K=\Cb$.
Let $f$ be a nonzero element of $\FF$. There exists $\chi'\in\Omega_\chi$
such that $f(\chi')\ne 0$ and the algebra $R_{\eps,\chi'}$ is isomorphic to
$R_{\eps,\chi}$.\\
{\bf Theorem 4.2}. Let $K=\Cb$, 
$R$ be a NQS-algebra
obeying Conditions CN1 and CN2,
and $\eps $ is an admissible specialisation of $C$.
Let $\pi$ be an irreducible representation with
central charakter $\chi$. Then  \\
1) $\dim(\pi)=l^{\frac{1}{2}\dim(\Omega_\chi)}$,\\
2) $\Omega_\chi$ is algebraic (i.e. it is Zariski open in its Zariski
closure),\\
3) the algebras $R_{\eps,\chi'}$ and $R_{\eps,\chi''}$
are isomorphic for any $\chi',\chi''\in \Omega_\chi$,\\
{\bf Proof}. 
For an irreducible representation $\pi$ with the central charakter 
$\chi$, we consider its kernel $I(\pi)$ in $R_\eps$.
This ideal is prime and the greatest $\DD$-stable ideal $I(\pi,\DD)$
in $I(\pi)$ is completely prime (see Theorem 3.3).
The ideal $m(\chi,\DD)$ coincides with $I(\pi,\DD)\bigcap Z_\eps$.

By Theorem 3.2, there exists $(\eps,\DD)$-admissible denominator set
$S:=S_\mu$ with empty intersection with
$I(\pi,\DD)$.
The ideal $I(\pi,\DD)$ admits localization
$\tID:= I(\pi,\DD)S_\eps^{-1}$ and 
$\tID\supset \DDJ$ (see Section 3). 
The subset $S_l:=\{h^l: h\in S\}$ is a denominator subset consisting of
$q$-commuting FA-elements and $\tR:=RS^{-1}=RS_l^{-l}$.
The subset $S_{l\eps}$ belongs to the center $Z_\eps$; it 
is a denominator subset in $R_\eps$ and 
$\tR_\eps:=R_\eps S_\eps^{-1}=R_{\eps}S_{l\eps}^{-1}$.
 The ideal $m(\chi,\DD)$ has empty entersection with $S_{l_\eps}$.
 We denote $\tZe=Z_\eps S_{l\eps}^{-1}$ and 
 $\tmD:= m(\chi,\DD)S_{l\eps}^{-1}=\tID\bigcap \tZ_\eps
 \supset \DDJ\bigcap\tZ_\eps$.
By Lemma 4.1, 
we may require 
$\chi(y_\eps^l)\ne 0$ for any $y\in S_l$.
Since $\pi$ is an irreducible representation, 
$\pi(y_\eps^l)=\chi(y_\eps^l)\cdot\rm{id}$.
The ideal $I(\pi)$ admits localization over $S_\eps$ and $\pi$ is an 
irreducible representation of $\tRe$

Recall $\tR_\eps:= \DDJ+Y_\eps$ (Section 3).
One can treat $\tR$ as a free left (and right) $Y$-module.
We form the free basis which consists of monomials 
(in the lexicographical order) of $\{\tx_i\}$. 
Denote by $\rho_S$ the natural projection
$\rho:\tR\to Y$.
The projection $\rho$ is a morphism of left (and right)
$Y$-modules. Similarly, $\rho_\eps:\tR_\eps\to Y_\eps$ is a
morphism of $Y_\eps$-modules
and $\rho_\eps\DD_{y^l}(a)=\DD_{y^l}\rho_\eps(a)$, for any $y\in G$ 
and $a\in \tR_\eps$. It follows that $\rho_\eps(W_\eps)=Z(Y)_\eps$.

The representation $\pi$ passes through $\rho_\eps$ and determined by
$$\nu_\alpha, 1\le \alpha\le 2r;\quad\lambda_\beta, 1\le \beta\le t;\quad
\xi_\gamma, 1\le \gamma\le p
$$
where $\pi(h_{\alpha\eps}^l)=\nu_\alpha\cdot\rm{id}$, 
$\pi(u_{\beta\eps})=\lambda_\beta\cdot\rm{id}$
and $\pi(z_{\gamma\eps})=\xi_\gamma\cdot\rm{id}$.

The ideal $\tI(\pi):=I(\pi)S^{-l}_\eps $ is the maximal ideal of 
$Z(Y_\eps)$ generated by all
$h_{\alpha\eps}^l-\nu_\alpha$, $u_{\beta\eps}-\lambda_\beta$, 
$z_\gamma-\xi_\gamma$.
We obtain $$\dim(\pi)=l^{2r}.\eqno(4.1)$$

Denote by $\lambda$ the character of $Z(Y)_\eps$ determined by $\pi$.
The ideal 
$${\frak m}_\pi:=\Ker(\lambda) =
\sum_{1\le \beta\le t}Y_\eps(u_{\beta\eps}-\lambda_\beta) + 
\sum_{1\le \gamma\le p}Y_\eps(z_\gamma-\xi_\gamma)\eqno (4.2)$$
is a maximal $\DD_Y$-stable ideal in $Y_\eps$. 

Character $\lambda$ obeys
the condition $\lambda\vert_{\DDJ\bigcap Z(Y)_\eps} = 0$.
We have
$$\tID\subset \DDJ+Y_\eps {\frak m}_\pi\subset \tI(\pi).$$
By Proposition 3.17, the middle ideal is $\DD$-stable.
It implies 
$$\tID=\DDJ+Y_\eps {\frak m}_\pi.\eqno (4.3)$$

Similarly, by (3.7), the central character $\chi$ also passes 
through $\rho_\eps$
and determined by $\nu_\alpha = \chi(h_{\alpha\eps}^l)$
and $\chi\vert_{W_\eps\bigcap\tZ_\eps}$.

We consider 
$${\frak M}_\chi:= \rm{Ker}\chi\vert_{W_\eps\bigcap\tZ_\eps}$$
and obtain
$$\tmD=\DDJ\bigcap\tZ_\eps +\Hle{\frak M}_\chi.\eqno(4.4)$$

Comparing (3.7) and (4.4),
we see that
the algebra $\tZe/\tmD$ is isomorphic (as a Poisson algebra)
to
$\Cb[h_{1\eps}^{\pm l},\ldots,h_{2t\eps}^{\pm l}]$
with the Poisson bracket
$$
\{h_{1\eps}^l,h_{2\eps}^l\}= m_1 l^2\eps^{-1}h_{1\eps}^lh_{2\eps}^l, \ldots
, \{h_{2r-1,\eps}^l,h_{2r\eps}^l\}= 
m_r l^2\eps^{-1}h_{2r-1,\eps}^lh_{2r\eps}^l.
$$ 
 The maximal spectrum of above Poisson algebra
has a single symplectic leaf.
It follows that the symplectic leaf $\Omega_\chi$ contains 
the subset $\OO$ which is Zariski-open in 
the Zariski closure $\MM_\chi:=\overline{\Omega}_\chi$.
It follows $\dim\Omega_\chi=\dim\OO=2r$ and, by (4.3), 
$\dim(\pi)=l^{\frac{1}{2}\dim(\Omega_\chi)}$. This proves 1).\\
2) For any $\chi'\in\MM_\chi-\Omega_\chi$, we see
$$
\dim\Omega_{\chi'}\le\dim(\MM_\chi-\OO)< 
\dim\OO=\dim\Omega_\chi.
$$
Then $\MM_\chi-\Omega_\chi=\{\chi'\in\MM_\chi:
 \dim\Omega_{\chi'}< 2r\}$.
 On the other hand, the subset $\MM_{<2r}:=\{\chi'\in\MM:
 \dim\Omega_{\chi'}<2r\}$
 is Zariski closed in $\MM$. The subset $\MM_\chi-\Omega_\chi$
 coincides with $\MM_{<2r}\bigcap \MM_\chi$ 
 and is Zariski-closed. Hence $\Omega_\chi$ is Zariski-open in $\MM_\chi$.\\
3) 
Consider the equivalence relation on $\Omega_\chi$ as follows
$\chi\cong\chi'$ iff $R_\eps/m(\chi)R_\eps\cong R_\eps/m(\chi')R_\eps$.
Any equivalence class $[\chi]=\{\chi':\chi'\cong\chi\}$
is an open subset of $\Omega_\chi$
in the topology of complex smooth  manifold [P3, Lemma 5.1].
The manifold $\Omega_\chi$ is connected.
This proves claim 3).$\Box$ \\ 
{\bf Theorem 4.3}.
Let $K$ be an algebraic closed field of zero 
characteristic.  Let $R$ and $\eps$ be  as in Theorem 4.2.
 Any two vertices $e_i,e_j$, $i\ne j$ of the 
 quiver of algebra $R_{\eps,\chi}$ are linked by the wedges $(e_i,e_j)$ and
 $(e_j,e_i)$.
In particular the quiver is connected.\\
{\bf Proof}.\\  
{\it Step 1}. Let us prove that all irreducible representations over a
common central character $\chi$ can be passed through suitable 
localization $\tR_\eps$ such that $\pi(\DDJ)=0$ for any $\pi$ over $\chi$.

For any irreducible representation $\pi$ there exists
$(\eps,\DD)$-admissible standard denominator subset $S:=S_\mu$ such that
$I(\pi,\DD)\bigcap S_{l\eps}=\emptyset$ and $\tI(\pi,\DD)\supset\DDJ$.
We may assume that $I(\pi)\bigcap S_{l\eps}=\emptyset$
(see Proof of Theorem 4.2).
As above $\chi$ is the central character of $\pi$.
The ideal $\tmD$ admits localization over $S_{l\eps}$
and $\tI(\pi)\supset R_\eps \tm(\chi)\supset R_\eps \tm(\chi,\DD)$.

For an other irreducible representation $\pi'$ over $\chi$,
we also have $I(\pi')\supset R_\eps m(\chi)\supset R_\eps m(\pi,\DD)$.
For the greatest $\DD$-stable ideal $I(\pi',\DD)$ in $I(\pi')$,  
we obtain $ I(\pi')\supset I(\pi',\DD)\supset R_\eps m(\chi,\DD)$
and $ I(\pi',\DD)\bigcap Z_\eps = m(\chi,\DD)$.
It implies that $I(\pi',\DD)\bigcap S_{l\eps}=\emptyset$ and
$I(\pi',\DD)\supset\{\tilde{x}_{i\eps}^l, i\in [1,\tilde{M}-\mu]\}$
(see Section 3). That is $I(\pi',\DD)$ contains $\DDJ$ and 
admits localization over $S_{l\eps}$. 
Ideal $I(\pi')$ also admits
localization over $S_{l\eps}$. 
This proves the claim of Step 1. \\
{\it Step 2.}
According to Step 1, any irreducible representation $\pi$ over
 $\chi$ is a representation of $\tR_\eps$
and its kernel contains $\DDJ$. 
Then 
$\DDJ\bmod m(\chi)$ is contained in the radical of 
$R_\eps/m(\chi)R_\eps$.

By (3.7),
$\pi$ lies over $\chi$ iff 
$\pi(h_{i\eps}^l)=\chi(h_{i\eps}^l)\cdot\rm{id}$ and
$$\pi\vert_{W_\eps\bigcap\tZ_\eps}=\chi\vert_{W_\eps\bigcap\tZ_\eps}\cdot
\rm{id}.\eqno(4.5)$$

Denote  $Z(Y_\eps)':=\rho_\eps(W_\eps\bigcap \tZ_\eps)$.
Since $\rho_\eps(W_\eps)=Z(Y)_\eps$, then
$Z(Y_\eps)'$ is a subalgebra of $Z(Y)_\eps$.

The character $\chi$ defines the character $\chi'$ on $Z(Y_\eps)'$ such 
that
$\chi'\rho_S(w_\eps)=\chi(w_\eps)$ for any $w_\eps\in W_\eps\bigcap \tZ_\eps$.
In particular, $\chi'(z_{i\eps})=\chi(z_{i\eps})$, and 
$\chi'(u_{i\eps}^l)=\chi(u_{i\eps}^l)$.

According to the proof of Theorem 4.2,
there exists 1-1 correpondence between irreducible 
representations over $\chi$
and characters $\lambda$ of $Z(Y)_\eps$ such that 
$$\lambda\vert_{Z(Y_\eps)'}=
\chi'\vert_{Z(Y_\eps)'}.\eqno(4.6)$$
We will say that such $\lambda$ is comparable with $\chi$.
In particular, $\lambda(z_{i\eps})=\chi'(z_{i\eps})$, and 
$\lambda_i^l=\lambda(u_{i\eps}^l)=\chi'(u_{i\eps}^l)=:\chi_i$.

Two charactors $\lambda,\lambda'$ over $\chi$ differ
$\lambda'_i=\eps_i\lambda_i$ ,$1\le i\le r$ where
 $\eps_1,\ldots,\eps_r$ are $l^{\rm{th}}$ roots
of unity.
Denote
$$
e_\lambda=l^{-r}\prod_{i=1}^r((\lambda_i^{-1}u_{i\eps})^{l-1}+
(\lambda_i^{-1}u_{i\eps})^{l-2} + \cdots + 1).
$$
The elements $\{e_\lambda\}$ obey $e_\lambda^2= e_\lambda$.
If $\lambda$ is comparable (resp. non-comparable) with $\chi$, then
$e_\lambda$ is a primitive idempotent corresponding $\pi$
(resp. is a zero element of $\tR_\eps/\tm(\chi)\tR_\eps$). 

By choice of $u_1,\ldots, u_t$ (see (3.5) and below), there 
exist $v_1,\ldots, v_t$ such that
$v_iu_j=q^{n_{ij}}u_jv_i$ where  $d:=\rm{det}(n_{ij})_{ij=1}^t\ne 0$
and $d$ is relatively prime with $l$ (see Definition 3.1).
For any system $(\eps_i,\ldots,\eps_r)$
of $l^{\rm{th}}$ roots of unity, there exists $v\in \tR_\eps$ 
such that 
$$vu_i=\eps_iu_iv.\eqno (4.7)$$
Let us prove that one can choose $v\notin \tm(\chi)\tR_\eps$.
Since the $\Ad$-action of the subgroup $U_\eps$ generated
by $u_{i\eps}$, $1\le i\le r$ is diagonalizable,
one can decompose 
$\tR_\eps = \tm(\chi)\tR_\eps\oplus V$ where
$V$ is some finite dimensional $\Ad_{U_\eps}$-stable subspace.
Consider the completion $\widehat{R_\eps}$ (resp. $\widehat{Z_\eps}$)
of $\tR_\eps$ (resp. $\tZ_\eps$) in the $\tm(\chi)$-adic
topology.
We have decomposition    
$\widehat{R_\eps}=\widehat{Z_\eps}\otimes \tR_\eps
\cong \widehat{Z_\eps}\oplus V $.
The $\Ad_{U_\eps}$-action
is identical in $\widehat{Z_\eps}$.
One can choose $v\in V$.

Put $\hat{v}:=v\bmod\tm(\chi)\tR_\eps$.
We have proved that $\hat{v}\ne 0$.
For $(\eps_1,\ldots,\eps_r)\ne (1,\ldots,1)$, the element  $v$ lies in
$\DDJ$.
The formula (4.7) implies that for different primitive 
idempotents $e_\lambda$,
$e_{\lambda'}$ of $\tR_\eps/\tm(\chi)\tR_\eps$
there exits an nonzero element $\hat{v}$ of the radical such that
$$\hat{v}e_\lambda=e_{\lambda'}\hat{v}.$$
 
The idempotents 
$\lambda$ and $\lambda'$ are linked by the wedge 
(as vertices of the quiver of algebra $\Rchi$) 
[Pie, 6.4]. 
$\Box$\\

\section{ On number of irreducible representations}

The goal of this section is to prove 
the statements on the number of irreducible representations over the common
central character.

We begin with the proof of the formula (5.1) 
for some ideal in iterated skew polynomial 
extension. The property (5.1) is well known for commutative rings
[AM, Corollary 10.18]. Notice that, in general, (5.1) 
is false for noncommutative 
iterated extensions (for instance, take $R=U({\frak g})$ 
for two-dimensional Lie algebra $[x,y]=y$ and $I=<x,y>$).\\
{\bf Lemma 5.1}
Let we have an iterated $q$-skew extension
$\Rg=\Rg_1\supset\ldots\supset \Rg_n=\Yg$ 
of the algebra $\Yg$ over the field $\Fg$;
$\Rg_i=\Rg_{i+1}[x_i;\tau_i,\delta_i]$ where $\tau_i$ is  a 
diagonal automorphism
of $\Rg_{i+1}$ and $\tau_i\delta_i=q_{i}\delta_i\tau_i$ 
with $q_{i}\in \Fg^*$. 
We impose the following requirements.\\
1) $\Yg$ is a free module 
over its center and a Noetherian domain;\\ 
2) any ideal of $\Yg$ is generated by its intersection 
with the center; \\
3) $\delta_i(\Yg)=0$;\\
4) any $\delta_i$ is locally nilpotent in $\Rg_{i+1}$. 
Let $\Ig$ be an ideal of $\Rg$ which contains $x_1,\ldots,x_n$.
Then 
$$\bigcap_{m=1}^\infty \Ig^m=0.\eqno (5.1)$$\\
{\bf Proof}. We shall prove by induction on $n$. If $ n=1$, then $\Rg=\Yg$.
Let $\{f_\alpha\}$ be the free basis of 
$\Yg$ over its center $Z(\Yg)$.
The ideal $\Ig$ is generated by its intersection with $Z(\Yg)$.
Then $\Ig=\{\sum c_\alpha f_\alpha: c_\alpha\in\Ig\bigcap Z(\Yg)\}$.
Hence
$\Ig^m\subset\{\sum c_\alpha f_\alpha: c_\alpha\in(\Ig\bigcap Z(\Yg))^m\}$.
The property (5.1) is true for $\Ig\bigcap Z(\Yg)$; it is true for $\Ig$.

Suppose that (5.1) is true for extensions of the length $\le n$.
Let us prove for an extension of the length $n+1$.
Let
$\Rg$ be the iterated extension of $\Yg$ that obeys the requirements
of the Lemma
$$\Rg=\Rg_*[x;\tau,\delta]\supset\Rg_*=
\Rg_1\supset\ldots\supset \Rg_n=\Yg.$$
By the induction hypothesis,
(5.1) is true for the ideal 
$\Ig_* = \Ig\bigcap\Rg_*$ of $\Rg_*$.

Since $x\in\Ig$, then $\delta(\Rg_*)\subset \Ig_*$.
Any element of $\Ig$ has the form
$r_0+xr_1+x^2r_2+\ldots $ where $r_0\in\Ig_*$ and $r_i\in\Rg_*$, $i\ge 1$.
Therefore,
$\Ig^m$ is the span of $x^kb_k$ ,
$$
b_k=\delta^{\alpha_1}(r_1)\cdots\delta^{\alpha_n}(r_n)j_1^{\beta_1}
\cdots j_t^{\beta_t}\eqno(5.2)
$$
where $k,\alpha_i,\beta_i,n,t$ are nonnegative integers, $r_i\in\Rg_*$, 
$j_i\in \Ig_*$
and 
  $$k+\alpha_1+\cdots+\alpha_n+\beta_1+\ldots+\beta_t\ge m.\eqno(5.3)$$
Suppose that $a\in\bigcap_{m=1}^\infty\Ig^m$ and $a\ne 0$.
Then $a=x^kb_k+x^{k+1}b_{k+1}+\ldots$, $b_k\ne 0$.
For any $m$ one can present $b_k$ in the form (5.2)
where
$\alpha_i $, $\beta_i$, $t$, $n$ depends on the choice of $m$
and (5.3) holds.

On the other hand, since $\bigcap\Ig_*^m=0$, there exists
$m_0$ such that 
$$b_k\in \Ig_*^{m_0}\qquad\mbox{and}\qquad b_k\notin\Ig_*^{m_0+1}.\eqno(5.4)$$
The condition (5.4) yields
$\beta_1+\ldots+\beta_t\le m_0$ and the number of nonzero $\alpha_i$ is also
$\le m_0$.

Recall that  $\delta $ is locally nilpotent $\tau$-derivation; 
there exists $N$ such that $\delta^N(x_i)=0$ for all $i$.
It implies $\delta^{nN}(\Rg_*)\subset\Ig_*^n$ for all $n$.
In particular, $\delta^{(m_0+1)N}(\Rg_*)\subset\Ig_*^{m_0+1}$.
 Since $b_k\notin\Ig_*^{m_0+1}$, then $\alpha_i<(m_0+1)N$ for any $i$.

 We conclude that the left side of inequality (5.3) is restricted as 
  $m$  tends to infinity.
This leads to a contradiction. The ideals $\Ig^m$ have zero intersection.
$\Box$.

Let $S:=S_\mu$ be the standerd denominator subset 
(see  Section 3).
Recall that after localisation of a NQS-algebra, we obtain a
an iterated $q$-skew extension 
$\tR:=
\tR_1\supset\cdots\supset
\tR_\tM\supset Y$
where $\tM:=M-k$ and $\tR_i\cong\tR_{i+1}[x;\ttau_i,\tdelta_i]$
(see (3.3)).
As above we denote by $\JJ=\JJ_S$ the lowest ideal of $\tR$ which contains
$\tx_i$ for all $i$.
Let $\{\QQ_1,\ldots,\QQ_m\}$ be the set of all minimal prime ideals 
over $\JJ$.
Denote 
$$X_1:=X_{1S}:=\{\QQ_i: \QQ_i\bigcap C=0\},$$ 
$$X_2:=X_{2S}:=\{\QQ_i: \QQ_i\bigcap C\ne 0\}.$$  
{\bf Proposition 5.2}. 
Suppose that $X_1\ne\emptyset$
and $\QQ\in X_1$.\\
1) The ideal $\QQ\bigcap Y$ is generated by $\QQ\bigcap \Zg$.\\
2) $\QQ= \JJ + \tR(\QQ\bigcap Y)$.\\
{\bf Proof}. 
The second statement is the easy corollary of the first.
Our goal is to prove statement 1).

First notice that any prime ideal of $X_1$ is completely prime.
[GL2, Theorem 3.2].
We will prove the statement by induction on $\tM$.
The statement if obviously true for $\tM=0$.
 
 We assume that the statement is true for an algebra of
 length $\tM$.  
Our aim is to prove the statement for 
$\tR$ of the length $\tM+1$.
Let $\tR_*$ be the subalgebra generated by $Y$ and all $\{\tx_i\}$ 
apart from the first, $\tR=\tR_*[\tx;\tilde{\tau},\tdelta]$, 
$\tilde{\tau}\tdelta=q^s\tdelta\tilde{\tau}$ with
$s\ne 0$ and $\JJ_*$ be the minimal ideal of $\tR_*$ which
contains $\{\tx_i\}$.

 Since completely prime ideal $\QQ\bigcap\tR_*$ of $\tR_*$ contains
 $\JJ_*$, then there exists some minimal prime
 ideal  $\QQ_*$ of $\tR_*$ such that $\QQ\bigcap\tR_*\supset\QQ_*\supset\JJ_*$.
 
Since $0=\QQ\bigcap C\supset\QQ_*\bigcap C$, then $\QQ_*\bigcap C=0$.
Whence $\QQ_*$ obeys the requirements of Proposition.
In partiqular, $\QQ_*$ is a completely prime ideal of $\tR_*$.  
The ideal $\qg:=\tR(\QQ_*\bigcap \Zg)$ of $\Zg$ is completely prime.
We  retain the former notations $\tR$, $Y$, $\QQ$
for
 $\frac{\tR}{\tR\qg}$, $\frac{Y}{Y\qg}$,$ \frac{\QQ}{\tR\qg}$.
By this agreement, $\QQ_*\bigcap Y=0$.
We obtain the natural projection 
$\pi_S:\tR_*\to Y$ with the kernal $\QQ_*$.
We denote by $B$ the denominator subset
$Y-\{0\}$.
The algebra $\Rg_*:=\tR_* B^{-1}$ 
is an iterated $q$-skew polynomial extension of
$\Yg:=YB^{-1}$. The ideal $\Ig_*:=\QQ_{*}B^{-1}$ 
of $\Rg_*$ obeys the requirements
of Lemma 5.1. One can choose the subset 
$\Psi\subset \{\tx_i\in\tR_*\}$ 
which form $\Yg$-basis of $\Ig_*$ over $\Ig_*^2$.
By (5.1), 
the set $\Psi^m$ which consists of products of arbitrary 
$m$ elements of $\Psi$,
generate $\Ig_*$ over $\Ig_*^{m+1}$. This implies that,
if an element of $Z(Y)$ 
commute with all elements of $\Psi$, 
then it lies in the center of $\tR$

We put 
$\tilde{\tau}(u_i)=q^{\alpha_i}u_i$
and denote by $\Zg_*$ the intersection of the center of $\tR_*$ with $Y$.\\
The subalgebra $\Zg_*$ is contained in $Z(Y)$. The following case take place.\\
Case 1.
$\Zg_*=\Zg$. There exist elements   
$\Phi:=\{v_{1},\ldots,v_{t}\}\subset\Psi$
such that $u_{i}v_{j}=q^{n_{ij}}v_{j}u_{i}$, 
$d:=\det(n_{ij})_{ij=1}^t\ne 0$.
We put $\Phi_*=\Phi$, if there is no $v_{j_0}\in\Phi$ such that
$u_iv_{j_0}=q^{\alpha_i}v_{j_0}u_i$  for all $i$. 
If the above $j_0$ exists, it is unique and we put $\Phi_*:=\Phi-\{v_{j_0}\}$.
\\
Case 2. $\Zg_*\ne\Zg$. 
One can suppose that $\Zg_*$ is generated by $\Zg$ and $u_t$.
 Remark $\alpha_t\ne 0$
(otherwise $u_t\in\Zg$). 
There exist elements
$\Phi_*=\{v_{1},\ldots,v_{t-1}\}\subset \Psi$
such that all $v_{i}$ commute with $u_{t}$, .
 $u_{i}v_{j}=q^{n_{ij}}v_{j}u_{i}$,
$d':=\det(n_{ij})_{ij=1}^{t-1}\ne 0$. 
In the Case 2, we put $\Phi:=\Phi_*\bigcup\{\tx\}$.\\
{\it Step 1}.
 We are going to prove that $\tdelta(v_{j})\in \QQ_*$
for any $v_j\in\Phi_*$. That is $b_j:=\pi_S(\tdelta(v_j))=0$
(for $v_j\in\Phi_*$).
 
Since $u_{i}v_{j}=q^{n_{ij}}v_{j}u_{i}$, then, using $\tdelta(u_i)=0$, we obtain
$
q^{\alpha_i}u_{i}\tdelta(v_{j})-q^{n_{ij}}\tdelta(v_{j})u_{i}= 0
$.
The element $u_i$ lies in the center $Z(Y)$ and is invertible.
We have $(q^{\alpha_i}-q^{n_{ij}})b_j = 0$.
Recall $v_j\in\Phi_*$; 
 there exists $i_0$ such that
$q^{\alpha_{i_0}}\ne q^{n_{i_0j}}$. It implies
$b_j=0$. This concludes Step 1.\\
{\it Step 2}. 
Recall that $\tilde{\tau}$ (but not any $\tilde{\tau}_i$) is an automorphism of $\tR$.
Then the ideal $\QQ$ (and $\QQ\bigcap Y$) is $\tilde{\tau}$-stable.
As for $\tilde{\tau}_i$, $1\le i$, this map is the automorphism
of $\tR_{i+1}$ (but not of $\tR$).
We are going to prove that
the ideal $\QQ\bigcap Y$ is $\tilde{\tau}_i$-stable for $v_i\in \Phi_*$.
It sufficies to verify that, for any generator $a\in\{\tx_i\}$ of $\tR_*$, 
the element
$b:=\pi_S(\tdelta(a))$ is $\tilde{\tau}_i$-eigenvector. 

Each $a$ is an $\Ad_G$-eigenvector,
then $b$ is also $\Ad_G$-eigenvector with the same eigenvalue.
Multiplying $a$ (and $b$) by suitable monomial 
$h_1^{m_1}\cdots h_{2r}^{m_{2r}}$, we may assume that $\Ad_gb=b$.
That is $b\in Z(Y)$.
 
Each generator $a$
is a FA-element (indeed $\FAq$-element) of $\tR$.
Take $v_i\in\Phi_*$. Since all generators are $\tilde{\tau}$-eigenvector, 
$\tilde{\tau}(v_i)=q^{\beta_i}v_i$.
There exists a polynomial
$f(t)= c_0t^N + c_1t^{N-1} +\cdots + c_N$, $c_0\ne 0$, $c_N\ne 0$,
$c_i\in C$ (which is 
decomplosable $f(t):=c_0\prod_{m=1}^N(t-q^{\gamma_m})$)
such that 
$$
c_0a^Nv_i+c_1a^{N-1}v_ia+\cdots+c_Nv_ia^N=0.\eqno(5.5)
$$
We act by  $\tdelta$  $N$ times on (5.5)  and  obtain
$$
 c'_0\tdelta(a)^Nv_i+c'_1\tdelta(a)^{N-1}v_i\tdelta(a)+
 \cdots + c'_Nv_i\tdelta(a)^N=
 0\bmod(\QQ^2_*)
 $$
 where $c'_j:=q^{\beta_i}c_j$. Then
$$
c'_0b^Nv_i+c'_1b^{N-1}v_ib+\cdots+c'_Nv_ib^N= 0\bmod(\QQ^2_*),
$$
$$
(c'_0b^N+c'_1b^{N-1}\tilde{\tau}_i(b)+\cdots+
c'_N\tilde{\tau}_i(b)^N)v_i= 0\bmod(\QQ^2_*).
$$
Each $v_i$ is an element of $\Yg$-basis of $\Ig$ over $\Ig^2$.
Therefore
$$
c'_0b^N+c'_1b^{N-1}\tilde{\tau}_i(b)+\cdots+c'_N\tilde{\tau}_i(b)^N=0,
$$
$$
\prod_{m=1}^N(b-q^{\gamma_m+\beta_i}\tilde{\tau}_i(b))=0.
$$

Since $Y$ is a domain, then $b$ is $\tilde{\tau}_i$-eigenvector.
The ideal $\QQ\bigcap \YY$ is $\tilde{\tau}_i$-stable.\\
{\it Step 3}.
Any ideal of $Y$ is generated by its intersection with the center;
 $\QQ\bigcap Y$ is generated by $\QQ\bigcap Z(Y)$.
 We have proved that
$\QQ\bigcap Y$ is $\tilde{\tau}_i$-stable for $v_i\in\Phi$, that is 
$\QQ\bigcap Y$ is generated by 
$\tilde{\tau}_i$-eigenvectors (for $v_i\in\Phi$).
All this eigenvectors have the form
$u_1^{m_1}\cdot u_t^{m_t} z$ where $z\in\Zg$. The elements $u_i$ are
 invertible;
the ideal $\QQ\bigcap Y$ is generated by $\QQ\bigcap \Zg$.
 $\Box$.\\
{\bf Notation 5.3}. For any $S:=S_\mu$ and any $\QQ\in X_{2S}$,
we consider the 
finite subset $\Eb_{S,\QQ}\subset K$ which consists of elements 
$\eps\in K$ such that $\QQ\bigcap C\vert_{q=\eps}=0$. 
We denote 
$$\Eb_S:=\bigcup_{\QQ\in X_{2S}} \Eb_{S,\QQ},$$
$$\Eb:=\bigcup_{S} \Eb_S.$$
Notice that the sets $\Eb_S$ and $\Eb$ is finite.

We consider specialisation of NQS-algebra $R_\eps$ 
at admissible root of unity.
Let $J:=J_S$ and
$\DDJ$ be the ideals of $R_\eps$ that were defined in Definition 3.1.
Let $P$ be some minimal prime ideal of $\tR_\eps$ over $J$.\\
{\bf Proposition 5.4}. $R$, $\eps$, $P$ as above.
Suppose that $\eps\notin \Eb_S$. Then \\
1) $P\bigcap Y_\eps$ is generated by $P\bigcap \Zg_\eps$,
$P=\DDJ + \tR_\eps(P\bigcap\Zg_\eps)$;\\
2) $P$ is a $\DD$-stable ideal.\\
{\bf Proof}.
The ideal $\PP:=\pi_\eps^{-1}(P)$ is prime and $\PP\bigcap C=(q-\eps)C$.
By definition of $P$, the ideal $\PP$ contains $\JJ$.
Then $\PP$ contains some minimal prime ideal $\QQ$ over $\JJ$.
If $\QQ\in X_{2S}$, then
$(q-\eps) C=\PP\bigcap C\supset\QQ\bigcap C$. Whence
$\QQ\bigcap C$ is zero at $q=\eps$. This leads to contradiction. 

Hence
$\QQ\in X_{1S}$.  
According to Proposition 5.2, $\QQ\bigcap Y$ is generated by $\QQ\bigcap \Zg$.
Then $\QQ=\JJ+\tR(\QQ\bigcap \Zg)$. Specialising modulo $q-\eps$, we
obtain 
 $\QQ_\eps = J+\tR_\eps\qg$
 where $\qg:=\QQ_\eps\bigcap\Zg_\eps$.
We have $P\supset\QQ_\eps\supset J$.
The ideal $\QQ_\eps$ is $\DD$-stable [P3, Lemma 3.12].
It  implies $\QQ_\eps\supset\DDJ$, $\QQ_\eps = \DDJ+\tR_\eps\qg$
and $\DDJ\bigcap Y_\eps\subset Y_\eps\qg$.

Recall that the ideal $P$ is prime; $P\bigcap\Zg_\eps$ is a prime 
ideal of $\Zg_\eps$.
There exists minimal prime ideal $\pg$ of $\Zg_\eps$ such that
$P\bigcap\Zg_\eps\supset\pg\supset\qg$.
We have
$P\supset\tR_\eps\pg$. Since $P\supset\QQ_\eps\supset\DDJ$,
then $P\supset\DDJ+\tR_\eps\pg\supset J$.
The middle ideal is 
prime, then $P=\DDJ+\tR_\eps\pg$. Similarly to Proposition 3.7,
$P$ is a $\DD$-stable ideal. $\Box$\\
{\bf Theorem 5.5}. Let $R$ be a NQS-algebra obeying Condetions CN1, CN2 
and $\eps$ be an admissible root of 1.
Suppose, in addition, that $\eps\notin\Eb$. Then the number of 
irreducible representations over central charector $\chi$ equals
to $l^t$ for some nonnegative integer $t$.
(To explicit the geometrical sense of $t$, see Theorem 5.7).\\
{\bf Proof}. As in the Section 4 we may assume that $\chi(S_{l\eps})\ne 0$.
The ideals $I(\pi)$ and $I(\pi,\DD)$ admit localization on $S_{l\eps}$.
After localization we obtain the ideals $\tI(\pi)$ and $\tI(\pi,\DD)$
of $\tR_\eps$ (see the proof of Theorem 4.2) and 
$\tI(\pi)\supset\tI(\pi,\DD)\supset \DDJ\supset J$.
For some minimal prime ideal $P$ (see Proposition 5.4):
$$\tI(\pi)\supset\tI(\pi,\DD)\supset P=
\DDJ+\tR_\eps\pg\supset J.$$
As above we put $\chi(z_{j\eps})=\xi_j$.
For maximal ideal 
$\Xi:=\sum_{i=1}^p\Zg_\eps (z_{j\eps}-\xi_j)$ of $\Zg_\eps$,
we see $\Xi\supset\pg\supset\DDJ\bigcap \Zg_\eps
\supset J\bigcap \Zg_\eps$.
 
Denote $\tR_\xi:=\frac{\tR_\eps}{\tR_\eps\Xi}$, 
$Y_\xi:=\frac{Y_\eps}{Y_\eps\Xi}\cong 
K[h_{1\eps}^{\pm 1},\ldots,h_{2r\eps}^{\pm 1},
u_{1\eps}^{\pm 1},\ldots,u_{t\eps}^{\pm 1}]$
and 
$J_\xi:=(\DDJ +\tR_\eps\Xi)\bmod \Xi = (J +\tR_\eps\Xi)\bmod\Xi$.
Notice that
$$\tR_\xi=J_\xi\oplus Y_\xi.\eqno (5.6)$$

The subset $B:=Y_\xi-\{0\}$ is a denominator 
subset in $\tR_\xi$. After localization we get
the algebras $\Rg_\xi:=\tR_\xi B^{-1}$, $\Yg_\xi:=Y_\xi B^{-1}$ and the ideal
$\Jg_\xi:=J_\xi B^{-1}$ which obeys the requirements of Lemma 5.1.
Therefore, $\bigcap_m \Jg_\xi^m =0$.
We retain the notations for generators $\tx_{i\eps}$, $h_{i\eps}$,
$u_{i\eps}$ of $\tR_\eps$ for there images in $\Rg_\xi$. 
Choose the basis $\Psi:=\{\tx_{i\eps}\}$ of $\Yg_\xi$-linear space 
$\Jg_\xi$ over
$\Jg_\xi^2$.
This basis generate $\Jg_\xi$ in the following sense:
 for any element $a$ of $\Jg$ and any $m\in\Nb$ there exists an expression $b$ 
 of the elements $\Psi$ with coefficients of $\Yg_\xi$ 
 such that $a=b\bmod\Jg^m$.
 A monomial of $\Yg_\xi$ lies in the center of $\Rg_\xi$
 if it commutes with all elements of $\Psi$.

The central character $\chi$ defines character
$\chi'$ of the subalgebra $Z(Y_\eps)'$ (see Section 4) which is contained 
in $Z(Y)_\eps$. Recall (see the proof of Proposion 4.4) that there exists 1-1 
correspondence between the irreducible representations  which lie over
$\chi$ and characters $\lambda$ of $Z(Y)_\eps$
such that $\lambda\vert_ {Z(Y_\eps)'}=\chi'\vert_{Z(Y_\eps)'}$.
After factorization over $\Xi$ we obtain the similar statement for
the subalgebra $Z(Y_\xi)':= Z(Y_\eps)'\Xi)$.

Let us show that 
$$Z(Y_\xi)'= K[u_{1\eps}^{\pm l},\ldots,u_{t\eps}^{\pm l}].\eqno (5.7)$$
Suppose that $y_0\in Z(Y_\xi)'$.
Then there exists an element $w\in (W_\eps\bigcap \tZ_\eps)\bmod \Xi$ such that
$w=j_0+y_0$ where $j_0\in J_\xi$ and $y_0\in Y_\xi$.

Let $a=\tx_{i_\eps}\in\Psi$.
Since $aw=wa$ and $aj_0-j_0a\bmod\Xi\in\Jg_\xi^2$, then
$aj_0-j_0a=y_0a-ay_0=(y_0-\tilde{\tau}_i(y_0))a$ modulo $\tR_\eps\Xi$ lies in
$\Jg_\xi^2$.
By definition of $\Psi$, $y_0=\tilde{\tau}_i(y_0)$.
Thus, $\Ad_{y_0}$  is identical on $\Psi$ and, therefore,  $y_0\in Z(R_\xi)$.
This proves (5.7).  According to (3.7),

$$Z(R_\xi)=(J_\xi\bigcap Z(R_\xi))\oplus
K[h_{1\eps}^{\pm l},\ldots,h_{2r\eps}^{\pm l},
u_{1\eps}^{\pm l},\ldots,u_{t\eps}^{\pm l}].
\eqno (5.8)$$
The number of irreducible representations over $\chi$ equales to $l^t$.$\Box$\\
{\bf Definition 5.6}. For a point $\chi\in\MM$, we denote
by $G(\chi)$ the Poisson subalgebra 
$\{a\in m(\chi):\{a,m(\chi)\}\subset m(\chi)\}$ in $Z_\eps$.
One can see
$G(\chi)\supset m(\chi)^2$. The finite dimensional Lie algebra
$\gog(\chi):= G(\chi)/ m(\chi)^2$ is called the stabilizer of $\chi$ [KM, 1.1].
\\
{\bf Theorem 5.7}. $R$, $\eps$, $t$  as in 5.5. 
For any central character $\chi$ the stabilizer 
$\gog(\chi)$ is a semidirect sum $\gog(\chi):=\jg+\tg$
where $\jg$ is an ideal and  $\tg$ is 
the toric subalgebra of dimension $t$.
\\
{\bf Proof}. 
We may suppose that the central character $\chi$ and all
$\pi$ over $\chi$ admit localization over $S_{l\eps}$ for some
standard denominator subset $S=S_\mu$.
After specialisation on $S_{l\eps}$, and, without loss of generality,
we factorize the algebra $\tR_\eps$
modulo the ideal $\Xi$ (see Theorem 5.5).

The cotangent subspace $T^*_\chi(\MM)$ is a span 
of the
images under $D:\tZ_\eps\to m(\chi)/m(\chi)^2$ of the
elements $h_{i\eps}^l$, $u^l_{i\eps}$, and $j_\eps \in 
J_\xi\bigcap \tZ_\xi$. 
We put $\jg:= D(J_\xi\bigcap Z(R_\xi))$ and
$\tg:= D(K[u_{1\eps}^{\pm l},\ldots,u_{t\eps}^{\pm l}])$.
The subalgebra  $\gog(\chi)$ decomposes into the direct 
sum of linear subspaces
 $\gog(\chi):=\jg+\tg$.
 
Since the ideal $J_\xi$ is $\DD$-stable (see Proposition 5.4),
the subspace $\jg$ is an ideal of $\gog(\chi)$.
Since $\Ad_{u_{i}}$ is diagonalizable, the subalgebra $\tg$ is toric.$\Box$\\
{\bf Conjecture 5.8}. The number of irreducible representations over $\chi$
is equal to $ l^{\rm{rank}\gog(\chi)}$.

\end{document}